\documentclass{birkjour}

\title{Open orbifold Gromov-Witten invariants of $[\mathbb{C}^3/\mathbb{Z}_n]$: localization and mirror symmetry}

\author{Andrea Brini}
\address{Section de Math\'ematiques, Universit\'e de Gen\`eve, 2-4 Rue du Li\`evre 1210, Gen\`eve, Switzerland}
\email{Andrea.Brini@unige.ch}
\author{Renzo Cavalieri}
\address{Department of Mathematics, Colorado State University, 101 Weber Building, Fort
Collins, CO 80523-1874}
\email{renzo@math.colostate.edu}

\usepackage{amsrefs}
\usepackage{amsfonts}
\usepackage{amssymb}
\usepackage{amscd}
\usepackage{fullpage}
\usepackage[labelformat=empty]{subfig}
\usepackage{booktabs}
\usepackage{xypic}
\BibSpec{article}{%
+{}{\sc \PrintAuthors} {author}
+{,}{ \textrm} {title}
+{,}{ \textit} {journal}
+{,}{ } {volume}
+{}{ \IfEmptyBibField{journal}{}{\PrintDate{date}}} {transition}
+{,}{ } {pages}
+{,}{ } {status}
+{.}{} {transition}
+{}{ Preprint\IfEmptyBibField{date}{}{ \PrintDate{date}}, available at \tt} {eprint}
+{.}{} {transition}
}

\BibSpec{book}{%
+{}{\sc \PrintAuthors} {author}
+{,}{ \textit} {title}
+{.}{ } {series}
+{,}{ } {volume}
+{.}{ } {publisher}
+{,}{ } {place}
+{,}{ } {date}
+{.}{} {transition}
}

\BibSpec{collection.article}{
+{} {\sc \PrintAuthors} {author}
+{,} { \textrm } {title}
+{,} { in \it \PrintConference} {conference}
+{,} { pp.~} {pages}
+{.} {\PrintBook} {book}
+{,} { \PrintDateB} {date}
+{.} {} {transition}
}

\BibSpec{innerbook}{
+{.} { \emph} {title}
+{.} { } {part}
+{:} { \emph} {subtitle}
+{.} { } {series}
+{,} { } {volume}
+{.} { Edited by \PrintNameList} {editor}
+{.} { Translated by \PrintNameList}{translator}
+{.} { \PrintContributions} {contribution}
+{.} { } {publisher}
+{.} { } {organization}
+{,} { } {address}
+{,} { \PrintEdition} {edition}
+{,} { \PrintDateB} {date}
+{.} { } {note}
}

\allowdisplaybreaks[3]

\newcommand{\Cstar}{{\mathbb{C}^\ast}}
\newcommand{\proj}{{\mathbb{P}^1}}

\newcommand{\smbzen}[2]{\overline{M}_{0,#1}(B\bZ_n, 0;#2)}
\newcommand{\smbzf}[2]{\overline{M}_{0,#1}(B\bZ_4, 0;#2)}

\renewcommand{\gcd}{\rm gcd}

\newcommand{\de}{{\partial}}
\newcommand{\hX}{{\widehat X}}

\newcommand{\bbN}{\mathbb{N}}
\newcommand{\bE}{\mathbb{E}}

\newcommand{\bbZ}{\mathbb{Z}}
\newcommand{\bZ}{\mathbb{Z}}
\newcommand{\bbR}{\mathbb{R}}
\newcommand{\bbC}{\mathbb{C}}
\newcommand{\bbP}{\mathbb{P}}
\newcommand{\bbF}{\mathbb{F}}
\newcommand{\bbQ}{\mathbb{Q}}
\newcommand{\bbT}{\mathbb{T}}

\def\bary{\begin{array}} 
\def\eary{\end{array}} 
\def\ben{\begin{enumerate}} 
\def\een{\end{enumerate}}
\def\bit{\begin{itemize}} 
\def\eit{\end{itemize}}
\def\nn{\nonumber} 

\newcommand{\cO}{\mathcal{O}}

\newcommand{\cC}{\mathcal{C}}
\newcommand{\DD}{\mathcal{D}}
\newcommand{\LL}{\mathcal{L}}
\newcommand{\cS}{\mathcal{S}}

\newcommand{\cW}{\mathcal{W}}
\newcommand{\cG}{\mathcal{G}}

\newcommand{\HH}{\mathcal{H}}

\newcommand{\cF}{\mathcal{F}}

\newcommand{\cM}{\mathcal M}

\newcommand{\OO}{\mathcal{O}}

\def\beq{\begin{equation}}                     %
\def\eeq{\end{equation}}                       %
\def\bea{\begin{eqnarray}}                     
\def\eea{\end{eqnarray}}

\def\IZ{{\mathbb Z}}
\def\IR{{\mathbb R}}
\def\IC{{\mathbb C}}
\def\IE{{\mathbb E}}

\theoremstyle{plain}

\newtheorem{che}{Check}
\newtheorem{thm}{Theorem}[section]
\newtheorem{lemma}[thm]{Lemma}
\newtheorem{prop}[thm]{Proposition}
\newtheorem{conj}[thm]{Conjecture}

\newtheorem*{conj*}{Conjecture}

\newtheorem*{cor*}{Corollary}

\theoremstyle{definition}

\newtheorem{rem}[thm]{Remark}
\newtheorem{defn}{Definition}
\newtheorem*{example}{Example}

\newcommand{\GIT}[1]{/\!\!/_{\kern-.2em #1 \kern0.1em}}

\renewcommand{\l}{\left}
\renewcommand{\r}{\right}

\begin{document}

\begin{abstract}

  We develop a mathematical framework for the computation of open orbifold Gromov-Witten invariants of $[\bbC^3/\bbZ_n]$, and provide extensive checks with predictions from open string mirror symmetry. To this aim we set up a computation of open string invariants in the spirit of Katz-Liu \cite{kl:open}, defining them by localization. The orbifold is viewed as an open chart of a global quotient of the resolved conifold, and the Lagrangian as the fixed locus of an appropriate anti-holomorphic involution. 
We consider two main applications of the formalism. After warming up with
the simpler example of $[\bbC^3/\bbZ_3]$, where we verify physical
predictions of Bouchard, Klemm, Mari\~no and Pasquetti \cite{Bouchard:2007ys, Bouchard:2008gu},
the main object of our study is the richer case of $[\bbC^3/\bbZ_4]$, where two different choices are allowed for the Lagrangian. For one choice, we make numerical checks to confirm the $B$-model predictions; for the other, we prove a mirror theorem for orbifold disc invariants, match a large number of annulus invariants, and give mirror symmetry predictions for open string invariants of genus $\leq 2$.
\end{abstract}

\maketitle

\section{Introduction}

In recent years, string-theoretic dualities have spurred a flurry of activity
in the Gromov-Witten theory of Calabi-Yau manifolds. 
In Physics,
Gromov-Witten theory comes naturally in two flavors: the closed topological $A$-model
gives rise to a (virtually) enumerative theory of compact Riemann Surfaces
mapping to the target space, whereas its open string counterpart, where the strings propagate with their boundary
constrained to certain submanifolds of the target (called $D$-branes), should
correspond to a mathematical counting problem of maps from a Riemann surface
with non-empty boundary. In the context of the open topological $A$-model on toric
Calabi-Yau threefolds \cite{Aganagic:2000gs, Aganagic:2001nx}, mirror 
symmetry techniques have been developed in the physics literature and
have led in recent times to a complete recursive formalism for the calculation
of ``open Gromov-Witten invariants''  \cite{Marino:2006hs, Eynard:2007kz,
  Bouchard:2007ys}.
This progress on the Physics side of the subject raised a host of new
mathematical challenges,
as the purported relationship of topological open string theory amplitudes with a
counting problem of maps from a Riemann surface with non-empty boundary posed
the problem of developing a suitable mathematical framework for the definition
\cite{Solomon:2006dx} and the effective calculation \cite{kl:open,
  Graber:2001dw} of  such invariants.  \\

This paper is concerned with the open Gromov Witten theory of a toric
Calabi-Yau orbifold (of dimension $3$). We develop a mathematical framework
for the computation of open orbifold invariants: Eq.~\eqref{locformula}
simultaneously defines and computes open invariants for any orbifold of the
form $[\mathbb{C}^3/\mathbb{Z}_n]$. The upshot of the formula is that open
invariants are controlled by degree $0$ closed $GW$ theory with descendants
and a combinatorial function, which we call  {\it disc function}, depending on the group action defining the orbifold.
Anyone familiar with Atyiah-Bott localization will immediately recognize that
in fact our formula can be readily extended to compute invariants for an
arbitrary toric orbifold, up to the usual localization combinatorial yoga. \\


We apply formula (\ref{locformula}) to confirm several predictions coming from mirror symmetry. In this area, physics-based predictions
have been even more sharply ahead of their $A$-model counterpart. In
particular, the combined effect of the relationship of topological open string
amplitudes with quasi-modular forms \cite{Aganagic:2006wq}, physical
expectations about their behaviour under variation of the K\"ahler structure
of the target manifold, and the recursive formalism of \cite{Bouchard:2007ys}
eventually led to a series of predictions for  open orbifold
Gromov-Witten invariants \cite{Bouchard:2007ys, Bouchard:2008gu,
  Brini:2008rh}. First we focus on the orbifold $[\mathbb{C}^3/\mathbb{Z}_3]$.

\begin{che}
Numerical computations for disc invariants for $[\mathbb{C}^3/\mathbb{Z}_3]$ confirm the mirror symmetry predictions of \cite{Bouchard:2007ys, Bouchard:2008gu}.
\end{che}

Our main case of study is however $[\mathbb{C}^3/\mathbb{Z}_4]$, where we have two different choices of Lagrangian: we call asymmetric the case in which the Lagrangian intersects one of the two axes that are quotiented effectively,  symmetric when it stems from the axis with nontrivial isotropy.
In the asymmetric case, we prove a mirror theorem for disc invariants.

\begin{thm}
The analytic part of the $B$-model asymmetric disc potential for
$[\mathbb{C}^3/\mathbb{Z}_4]$ at framing 1 coincides, for positive winding numbers and up to signs, with the generating function of orbifold Gromov-Witten disc invariants obtained from (\ref{locformula}).
\end{thm}

For annulus invariants we are only able to perform numerical checks.

\begin{che}
Numerical computations for asymmetric annulus invariants for $[\mathbb{C}^3/\mathbb{Z}_4]$ agree with the mirror symmetry predictions.
\end{che}
This check is particularly interesting as the mirror symmetry computations involves non-trivially the fact that the $B$-model annulus potential is a quasi-modular form of $\Gamma(2) \subset SL(2, \bbZ)$. We might then regard this as an {\it a posteriori} $A$-model check of the relationship between generating functions of Gromov-Witten invariants of Calabi-Yau threefolds and modular forms in the context of {\it open} invariants.
In the asymmetric case, we also give mirror symmetry predictions for open
invariants in genus $\leq 2$. \\

\begin{che}
Numerical computations for symmetric disc invariants agree with the mirror
symmetry predictions. 
\end{che}

\subsection{Plan of the paper}
The paper is organized as follows. In Section \ref{sec:katzliu} we set up the
computation of open orbifold invariants of $[\bbC^3/\bbZ_n]$ by viewing the
orbifold as an open chart of a global quotient of the resolved conifold
$\cO_{\bbP^1}(-1) \oplus \cO_{\bbP^1}(-1)$, and the Lagrangian as the fixed
locus of an appropriate anti-holomorphic involution. Section
\ref{sec:Bmodel} aims at a self-contained review of the $B$-model setup of \cite{Aganagic:2000gs,
  Aganagic:2001nx, Aganagic:2006wq, Bouchard:2007ys, Eynard:2007kz} for a mathematical
audience and prepares the ground for the mirror symmetry computations in the
rest of the paper. In Section \ref{sec:c3z3} we consider the case of
$[\bbC^3/\bbZ_3]$ which was considered from the $B$-model point of view by
Bouchard, Klemm, Mari\~no and Pasquetti in \cite{Bouchard:2007ys,
  Bouchard:2008gu}; we move in Section \ref{sec:z4} to our main case of study: $[\bbC^3/\bbZ_4]$.  Finally, we collect in the Appendix a few technical results about the Eynard-Orantin recursion and its relationship with quasi-modular forms, and list part of the results of our $B$-model computations of higher genus open string potentials for $[\bbC^3/\bbZ_4]$.

\subsection*{Acknowledgements} 
We would like to thank G.~Bonelli, V.~Bouchard, P.~Johnson, M.~Mari\~no, S.~Pasquetti, Y.~Ruan and A.~Tanzini for 
discussions. We want to thank in particular M.~Manabe for correspondence after the
appearance of his manuscript \cite{Manabe:2009sf}, whose results partially overlap with
those of Section \ref{sec:c3z4basymm}, and we are grateful to him for
kindly acknowledging our $B$-model computations of open orbifold invariants
prior to publication. We would also like to thank AIM  and the organizers of the workshop ``Recursion structures in topological string
theory and enumerative geometry'' in Palo Alto (June 2009) where part of this work was carried out. 
The first author is supported by a post-doc fellowship of the Fonds National Suisse
 (FNS); partial support from the
European Science Foundation Programme ``Methods of Integrable Systems,
Geometry, Applied Mathematics'' (MISGAM)
and the
``Progetto Giovani 2009'' grant ``Teoria di stringa topologica e sistemi
 integrabili'' of the Gruppo Nazionale per la Fisica Matematica (GNFM-INdAM)
 is also acknowledged.

\section{The $A$-model side}
\label{sec:Amodel}

\subsection{Open Gromov-Witten invariants following Katz and Liu}
\label{sec:katzliu}

In \cite{kl:open}, Katz and Liu propose a tangent/obstruction theory for the moduli
space of open stable maps which parallels the construction in ordinary
Gromov-Witten theory. Consider an almost K\"ahler manifold $(X,J,\omega)$,
 a Lagrangian $L\subset X$, a class $\beta \in H_2(X,L)$ and classes $\gamma_i \in H_1(L)$ such that $\sum \gamma_i=\partial \beta$.  The sheaves of the obstruction theory (here described in terms of their fiber over a smooth moduli point $(\Sigma,f)$ of $\overline{M}_{g,h}(X,L|\beta;\gamma_1,\ldots,\gamma_h)$) fit in the exact sequence:
\begin{align}\label{ot}
0\to H^0(\Sigma,\partial \Sigma,T_\Sigma, T_{\partial\Sigma}) \to
H^0(X,L,f^\ast T_X, (f_{|\partial\Sigma})^\ast T_L) \to \mathcal{T}^1 \to \nonumber \\
 H^1(\Sigma,\partial \Sigma,T_\Sigma, T_{\partial\Sigma}) \to
H^1(X,L,f^\ast T_X, (f_{|\partial\Sigma})^\ast T_L) \to \mathcal{T}^2 \to 0
\end{align}
The expected dimension is:
\begin{equation}
\mbox{rk}\ \mathcal{T}^1-\mbox{rk}\ \mathcal{T}^2 = \mu(f^\ast T_X, (f_{|\partial\Sigma})^\ast T_L)-(\dim X-3)\chi(\Sigma),
\end{equation}
where $\mu$ denotes the generalized Maslov index of the real sub-bundle $(f_{|\partial\Sigma})^\ast T_L \subset f^\ast T_X$ \cite[Section 3.7]{kl:open}.
In the case that $X$ is a complex manifold and $L$ is the fixed locus of an
anti-holomorphic involution, the complex double of $(f^\ast T_X,
(f_{|\partial\Sigma})^\ast T_L)$ is $f^\ast_\mathbb{C} T_X$ and the Maslov
index coincides with the first Chern class of the latter bundle. Hence, for
$X$ a Calabi-Yau threefold, we obtain a moduli space of virtual dimension
$0$. With the additional assumption that the moduli space is endowed with a well
behaved torus action, Katz and Liu propose the existence of a virtual cycle,
and give an explicit formula for its localization to the fixed loci of the
torus action. Such cycle {\it does} depend on the torus action: different choices of
action lead to different enumerative invariants. \\

Next, Katz and Liu specialize to $X$ the resolved conifold, that is, the total space of
$\OO_{\proj}(-1)\oplus \OO_\proj(-1)$, and the Lagrangian $L$ being the
fixed locus of the anti-holomorphic involution $A: (z,u,v)\mapsto (1/\bar{z},
\bar{z}\bar{v},\bar{z}\bar{u})$, where we use local coordinates ($z$, $u$,
$v$) for a chart around $0\in \proj$. The standard circle action on  the base
$\proj$ preserves the equator ($L\cap \proj$). An extension of the circle
action to a $\Cstar$ action and lifting of the torus action to the total space
of the resolved conifold is compatible with the antiholomorphic involution if
it has Calabi-Yau weights. Any such choice, say with weights  $(\hbar, -a\hbar, (a-1)\hbar)$ ($\in H^\ast_{\mathbb{C}^\ast}(pt.)\cong H^\ast_{S^1}(pt.)=\mathbb{C}[\hbar]$)
over $0$, determines uniquely a real line bundle inside
$(T_L)|_{\mathrm{equator}}$, and this topological data in turn determines the
virtual cycle used to compute open invariants. The fact that the invariants
are not intrinsic to the geometry of ($X$, $L$) matches the physical
expectations from large $N$ duality \cite{Gopakumar:1998ki}, and in
particular makes $a$ the natural closed-string counterpart of the framing ambiguity of knot invariants in
Chern-Simons theory \cite{Witten:1988hf}. \\

The torus action on the target induces a torus action on the moduli space of
maps, and the fixed loci are easy to understand. The restriction of the
virtual cycle to the fixed loci is evaluated using sequence (\ref{ot}). Before
describing these steps in detail, we
recall the properties of two bundles that play a special role in the
restriction of the virtual cycle to the fixed loci.

\subsubsection{The bundles $L(2m)$, $N(m)$} 
\label{LmNm}

We describe two Riemann-Hilbert bundles on $(D^2,S^1)$ that play a special role in our story. 
For $m>0$ consider the bundle $\OO_\proj(2m)$ and the anti-holomorphic involution $\sigma:(z,u)\mapsto (1/\bar{z}, -\bar{z}^{-2m}\bar{u})$. The fixed locus for $\sigma$ is a real sub-bundle of the restriction of $\OO_\proj(2m)$ to the equator. We abbreviate Katz and Liu and define:
\begin{equation}
\label{lm}
L(2m) = (L(2m), L(2m)_\bbR) := (\OO_\proj(2m)_{|D^2}, \OO_\proj(2m)_{|S^1}^\sigma).
\end{equation}
The global sections of $L(2m)$ are by definition the $\sigma$-invariant sections of $\OO_\proj(2m)$, and they can be embedded torus equivariantly into the sections of the complex bundle $\OO_\proj(m)$:
\begin{eqnarray}
\label{idlm}
H^0(L(2m)) &\hookrightarrow &  H^0(\OO_\proj(m))\nonumber \\
\sum_{j=0}^{m-1} (a_j z^j -\bar{a}_{j}z^{2m-j}) +ibz^m &\mapsto& \sum_{j=0}^{m-1} a_j z^j+bz^m,
\end{eqnarray} 
with $a_j\in \bbC, b\in \bbR$. Therefore the weights of the torus action for the left hand side can be computed in terms of the weights for the right hand side.

\begin{rem}
The identification (\ref{idlm}) chooses an orientation for the space of global sections of $L(2m)$. 
\end{rem}
\begin{rem}
There is an abuse of notation in saying ``torus equivariantly'', since a real torus acts on the left vector space, while the complex torus $\Cstar$ acts on the right. Here we identify the circle with $U(1)\subset \Cstar$. For $w\in\bbZ$, we identify the real $2$ dimensional $S^1$-representation corresponding to rotation by $w\theta$ with the one dimensional complex $\Cstar$ representation corresponding to multiplication by $\alpha^w$, and we give both weight $w\hbar$. By weight $0$ we mean the trivial representation, which is one real dimensional in the real case, and one complex dimensional in the complex case.   
\end{rem}

For $m>0$, now consider  $\OO_\proj(-m)\oplus \OO_\proj(-m)$. The anti-holomorphic involution $\sigma:(z,u,v)\mapsto(1/\bar{z},\bar{z}^m\bar{v},\bar{z}^m\bar{u})$ fixes a two dimensional real sub-bundle on the equator that we use to define $N(m)$:
\beq
\begin{aligned}
\label{nm}
N(m) = (N( m), N( m)_\bbR) := ((\OO_\proj(-m)\oplus\OO_\proj(-m))_{|D^2}, \\ (\OO_\proj(-m)\oplus\OO_\proj(-m))_{|S^1}^\sigma).
\end{aligned}
\eeq
The  sections of the first cohomology group of $N(m)$ are by definition the $\sigma$-invariant sections of $H^1(\OO_\proj(-m)\oplus \OO_\proj(-m))$, and an orientation is chosen by the torus equivariant identification with the sections of $H^1(\OO_\proj(-m))$:
\begin{eqnarray}
\label{idnm}
H^1(N(m)) &\to &  H^1(\OO_\proj(-m))\nonumber \\
\left(\sum_{j=1}^{m-1} \frac{\bar{a}_j}{z^{d-j}},\sum_{j=1}^{m-1} \frac{a_j}{z^j} \right)  &\mapsto& \sum_{j=1}^{m-1} \frac{a_j}{z^j}.
\end{eqnarray} 

\subsection{The orbifolds $[\IC^3/\IZ_n]$}
\label{sec:c3zn}

\subsubsection{The geometric set-up}
In this section we specialize the framework of \cite{kl:open} to the case of the ``orbifold vertex", deriving general formulas for open Gromov-Witten invariants in terms of the closed full descendant Gromov-Witten potential. Identify $\IZ_n$ with the multiplicative group of $n$-th roots of unity, set $\varepsilon= \mbox{e}^{\frac{2\pi i}{n}}$ and consider the quotient 
$\mathfrak{X}=[\IC^3/\IZ_n]$ by a Gorenstein action:
\begin{equation}
\label{goraction}
\varepsilon \cdot (x_0, x_1, x_2 ) = (\varepsilon^{\alpha_0} x_0, \varepsilon^{\alpha_1} x_1, \varepsilon^{\alpha_2} x_2 ),
\end{equation}
with $\alpha_0+\alpha_1+\alpha_2 \equiv 0$ (mod $n$). We wish to view our orbifold as an open chart of a global quotient of the resolved conifold:
$$
\mathfrak{X} \subset [\mathcal{O}(-1)\oplus \mathcal{O}(-1)/\IZ_n]
$$
Recall that $\mathcal{O}_{\bbP^1}(-1)\oplus \mathcal{O}_{\bbP^1}(-1)$ can be given local coordinates $(z,u,v)$ at $0$, $(z',u',v')$ at $\infty$ and the transition functions are: $z'=1/z, u'=uz,v'=vz$. Making the identification $(x_0,x_1,x_2)=(z,u,v)$, 
the action (\ref{goraction}) on the chart centered at $0$ induces an action on  the chart at $\infty$:
$$
\varepsilon\cdot(z',u',v'):= (\varepsilon^{-\alpha_0} z',\varepsilon^{\alpha_0+\alpha_1} u', \varepsilon^{\alpha_0+\alpha_2} v').
$$
Define an anti-holomorphic involution:
$$
\sigma(z,u,v)= (1/\bar{z},\bar{z}\bar{v},\bar{z}\bar{u}).
$$
The fixed locus of $\sigma$ is a Lagrangian  $\LL$ with topology $S^1\times \IR^2$ and explicit equation:
$$
\LL=\{(e^{i\theta}, u, \bar{z}\bar{u})\}
$$
One checks that $\varepsilon \cdot\sigma(-)=\sigma(\varepsilon \cdot -)$: hence $\sigma$ descends to the quotient defining a Lagrangian $\LL\subset \mathfrak{X}$.
We want a $\Cstar$ action on the total space of the resolved
conifold, which lifts the canonical action on $\proj$, descends to the
quotient and is compatible with the anti-holomorphic involution (that is, $U(1) \subset \Cstar$ preserves the Lagrangian): 
$$\sigma(\lambda P)=1/\bar{\lambda} \sigma(P).$$

Any Calabi-Yau action, i.e. an action where the sum of the three weights for the tangent space of a fixed point equals zero, satisfies these requirements. Since the weight $s_0$ is canonically linearized via the standard action on the tangent bundle to the doubled (orbi)-disc, the weights are determined up to the choice of a free parameter. It is convenient to use fractional weights for the induced action on the quotient:
$$
(s_0,s_1,s_2)=\left(\frac{\hbar}{n_{\mathrm{eff}}}, -{a\hbar}, { a\hbar }-\frac{\hbar}{n_{\mathrm{eff}}}\right),
$$
 The parameter $a$ should then correspond (up to an ``integer/$n$''
 translation) to the large $N$ dual of the {\it framing} ambiguity of
 Chern-Simons knot invariants \cite{Witten:1988hf, Aganagic:2001nx}. To keep the notations lighter in the general formulas we continue to use $(s_0,s_1,s_2)$, implicitly intending them as functions of the framing $a$ as in the above equation.

  
\subsubsection{The fixed loci}

The fixed maps for the torus action consist of a compact curve, possibly with twisted marks, with a collection of orbi-discs attached, depicted in Figure \ref{fig:fix}. The origin of the discs can be twisted, and the corresponding attaching point on the compact curve is twisted by the opposite character. The compact curve contracts to the (image of the) origin, and the discs are mapped rigidly to the zero section of   $[\mathcal{O}(-1)\oplus \mathcal{O}(-1)/\IZ_n]$, with their boundary wrapping around the equator. We describe such a map $f$ via the universal diagram of its complex doubling. Doubling the disc we obtain an  orbi-sphere $\mathcal{C}$ with $0$ a $k$-twisted point, $\infty$ a $(-k)$-twisted point. A fixed map $f$ of degree $d$ is then described by the following diagram: 
\begin{equation}
\label{univdiag}
\xymatrix{
\bigcup_{j=1}^{t_{\mathrm{in}}}(\proj,x_j) \ar[dd]_{X=x_j^{t_{\mathrm{eff}}}}\ar[rrr]^{F_\IC}& & & (\proj,z)\ar[dd]^{Z=z^{n_eff}} \\
   & \cC \ar[r]^{f_\IC} \ar[dl]& [\proj/\IZ_n] \ar[dr]& \\
   (\proj, X) \ar[rrr]^{Z=X^d}& & & (\proj, Z) \\ 
}
\end{equation}
where
\begin{align*}
n_{\mathrm{in}}= \gcd(\alpha_0,n) & & n_{\mathrm{eff}}= n/\gcd(\alpha_0,n)\\
t_{\mathrm{in}}= \gcd(k,n) & & t_{\mathrm{eff}}= n/\gcd(k,n)\\
& F_\IC: \{z=
x_j^{dn_{\mathrm{in}}/t_{\mathrm{in}}}\} &
\end{align*}
and the diagonal maps are the projections to the coarse moduli spaces. The $\IZ_n$ action on the upper-left collection of $\proj$'s is defined as follows: if $p_j\in (\proj,x_j)$,then  $\varepsilon \cdot p_j=p_{j+1} \in (\proj,x_{j+1})$ and
$$ 
x_{j+1}(p_{j+1})= \varepsilon x_j(p_j).
$$

\begin{rem}

\begin{figure}
	\centering
		\includegraphics[width=0.75\textwidth]{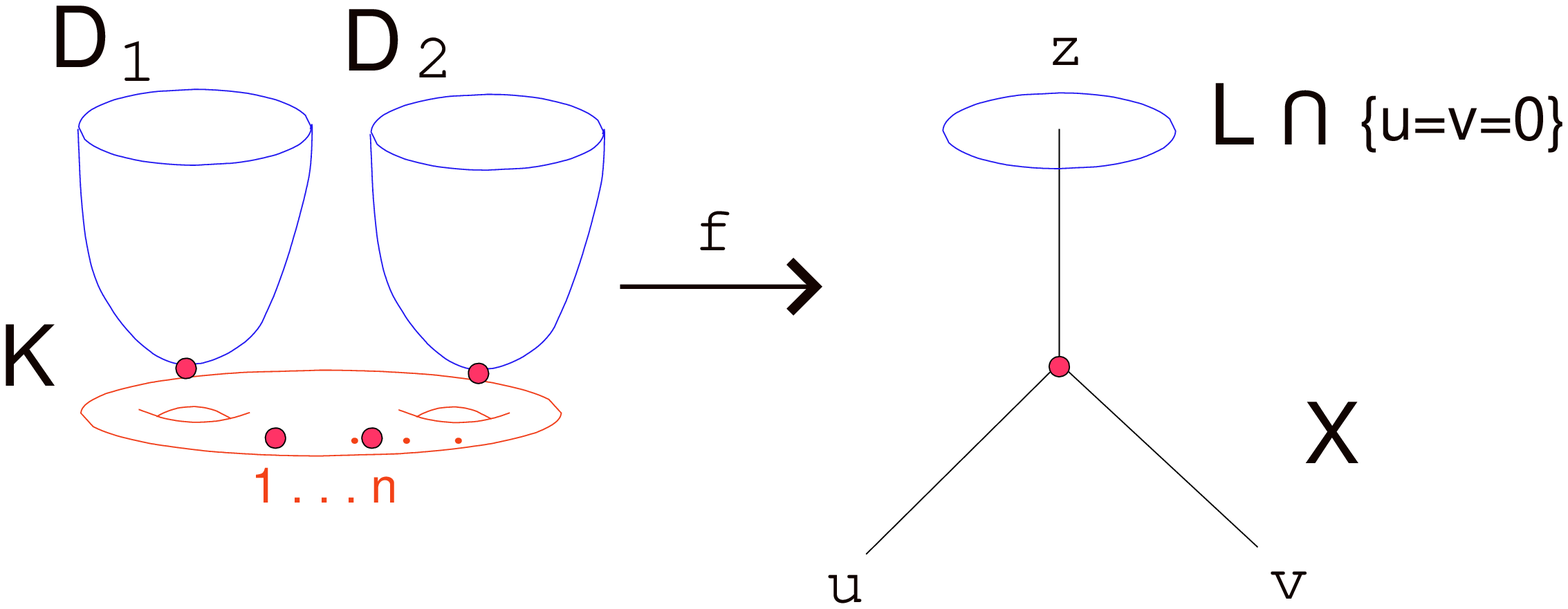}
	\caption{A fixed map: the compact (orbi)-curve $\mathcal{K}$ , hosting all marks, contracts to the ``origin'' of  $\mathfrak{X}$. Two discs, attached to $\mathcal{K}$ at possibly stacky points, map to the ``$z$-axis" their boundaries winding around the intersection of the Lagrangian $\mathcal{L}$.}
	\label{fig:fix}
\end{figure}

It is immediate to check that the above diagram is $k$-twisted equivariant ($F_{\IC}(\varepsilon^{t_{\mathrm{in}}} \cdot x_j)=\varepsilon^k \cdot F_{\IC}(x_j)$) if and only if:
\begin{equation}
\label{degtwist}
d \equiv k \frac{\alpha_0}{n_{\mathrm{in}}} \ \ \ (\mbox{mod}\ n_{\mathrm{eff}}),
\end{equation}
and therefore this numerical condition between degree and twisting must hold
for  $f$ to exist. Note also that equation (\ref{degtwist}) guarantees that
the degree of $F_\mathbb{C}$ is always integer.
\label{rmk:Amonodr}
\end{rem}

\subsubsection{The obstruction theory}

In this section we give a formula for the restriction of the obstruction theory (\ref{ot}) to a particular fixed locus in terms of the combinatorial data of the fixed locus. We give a careful treatment of the disc contribution, since that is essentially the part which is new. We denote by $d$ the winding degree, and by $k$ the twisting at the center of the disc.  
\begin{itemize}
\item {\bf compact curve}: the contribution by a contracting compact curve is given by the equivariant euler class of three copies of the dual of the appropriate $\alpha_i$-character sub-bundles of the Hodge bundle, linearized with the weights of the torus action. Notation and further explanation can be found, for example, in \cite[Section 2.1]{cc:c3z3}:
\begin{equation}
e^{\mathrm{eq}}(\IE^\vee_{\alpha_0}(s_0)\oplus\IE^\vee_{\alpha_1}(s_1)\oplus\IE^\vee_{\alpha_2}(s_2)).
\end{equation} 
\item {\bf node}: each node contributes a torus weight for any direction that the twisting makes invariant (i.e. $s_i$ if $k\alpha_i \equiv 0$ (mod $n$)), a denominator corresponding to smoothing the node. There is a gluing factor of $n$ (carefully discussed in \cite[Section 1.4]{cc:c3z3}). And finally we include an automorphism factor at the denominator to cancel the automorphisms of the disc. Define $\delta_i:=s_i^{\delta^{[n]}_{ka_i,n}}= s_i$ if $k\alpha_i \equiv 0$ (mod $n$), and $1$ otherwise. Then the contribution is: 
\begin{equation}
\delta_{0}\delta_{1}\delta_{2}\cdot\frac{1}{\frac{\hbar}{d}- \psi }\cdot n \cdot \frac{1}{\frac{\hbar}{d}}
\end{equation}  
\item {\bf disc}: a degree $d$, $k$-twisted at $0$,  fixed map $f$ from a disc has $n_{in}$ automorphisms, an $\hbar/d$ factor for infinitesimal automorphisms, and a contibution from the pull-back of the tangent bundle to $\mathfrak{X}$:
\begin{equation}
\frac{1}{n_{in}}\frac{\hbar}{d} \cdot e^{\mathrm{eq}}(R^\bullet_\ast f^\ast(T_\mathfrak{X},T_\LL)):=D_k(d,a). 
\end{equation} 
\end{itemize}

In the remainder of this section we give an explicit discussion and derive
formulas for this last contribution, which we dub {\bf disc function}. We study $f^\ast(T_\mathfrak{X},T_\LL)$ by studying its pull-back to the universal diagram. Splitting the bundle   into its tangent and normal component to the $x_0$ direction  
we have:
\begin{equation}
F^\ast(T_{\IC^3},T_L) = \bigoplus_{j=1}^{t_{\mathrm{in}}} L\left(2\frac{dn_{\mathrm{in}}}{t_{\mathrm{in}}}\right) \oplus \bigoplus_{j=1}^{t_{\mathrm{in}}}N\left(\frac{dn_{\mathrm{in}}}{t_{\mathrm{in}}}\right).
\end{equation}

This bundle (over an irreducible component of the fixed locus) is trivial but not equivariantly trivial. Its weights are computed via the identification discussed in Section \ref{LmNm}, as in \cite{kl:open}.
We must take this process one step further and select the sections that descend to the orbifold bundles, i.e. that are invariant under the $\IZ_n$ action. Referring to diagram (\ref{univdiag}) to identify the appropriate local coordinates, and defining $\mathbf{{x}^{\ell}}= \sum_{j=1}^{t_{\mathrm{in}}} x_j^\ell$, 
 we have:

\bea
\label{tginvsec}
H^0(L(2m))^{\IZ_n} & \stackrel{T-equiv}{\cong}
&\bigoplus_{j=1}^{t_{\mathrm{in}}}\left\langle \frac{\partial}{\partial
  z},{{{x_j}}}\frac{\partial}{\partial
  z},{{{x_j}^{2}}}\frac{\partial}{\partial z},\ldots,
{{{x_j}^{m-1}}}\frac{\partial}{\partial z}\right\rangle^{\IZ_n} \nn \\ &=&
\left\langle\left.{\mathbf{{x}^{\ell}}}\frac{\partial}{\partial
  z}\right|0\leq\ell<m,\ell \equiv m\  (\mbox{mod}
\ t_{\mathrm{eff}})\right\rangle. \nn \\
\eea

\begin{rem}
The section ${\mathbf{x}^{m}}\frac{\partial}{\partial z}$, corresponding to the pull-back of $z\frac{\partial}{\partial z}$, does not appear in the above list as it is acted upon trivially both by the torus and $\IZ_n$. This also explains the congruence in the last equality of (\ref{tginvsec}).
\end{rem}

For the normal part of the obstruction theory:

\bea
H^1(N(m))^{\IZ_n} & \stackrel{T-equiv}{\cong} & \bigoplus_{j=1}^{t_{\mathrm{in}}}\left\langle {{\frac{1}{x_j}}},{\frac{1}{{x_j}^{2}}},\ldots, {\frac{1}{{x_j}^{m-1}}}
\right\rangle^{\IZ_n} \nn \\ &=&
\left\langle{\mathbf{{x}^{\ell}}}\left|-m<\ell<0,\ell \equiv
\frac{k\alpha_1}{t_{\mathrm{in}}}\  (\mbox{mod}
\ t_{\mathrm{eff}})\right.\right\rangle. \nn \\
\eea
 
To compute the torus weights of the invariant sections, we look at the weights over $0$:
\begin{enumerate}
	\item The section $\frac{\partial}{\partial z}$ has weight $s_0=\hbar/n_{\mathrm{eff}}$.
	\item The section $x_j$ has weight $-\frac{t_{\mathrm{in}}}{nd}\hbar$.
	\item The trivializing section for the bundle $N(m)$ has weight $s_1$.
\end{enumerate}
Piecing everything together, we obtain:
\begin{equation}
e^{\mathrm{eq}}\left(H^0(L(2m))^{\IZ_n}\right)=\prod_{r=1}^{\lfloor \frac{d}{n_{\mathrm{eff}}}\rfloor} \frac{\hbar r}{d}=\lfloor \frac{d}{n_{\mathrm{eff}}}\rfloor!\left(\frac{\hbar}{d}\right)^{\lfloor \frac{d}{n_{\mathrm{eff}}}\rfloor}
\end{equation}

\begin{align}
e^{\mathrm{eq}}\left(H^1(N(m))^{\IZ_n}\right)=\prod_{r=1}^{\lfloor \frac{d}{n_{\mathrm{eff}}}-\frac{1}{t_{\mathrm{eff}}}+\left\langle\frac{k\alpha_1}{n}\right\rangle\rfloor}\left( s_1 -\frac{\hbar}{d}\left\langle\frac{k\alpha_1}{n}\right\rangle+\frac{\hbar r}{d}\right)= \nonumber \\ \left(\frac{\hbar}{d}\right)^{\lfloor \frac{d}{n_{\mathrm{eff}}}-\frac{1}{t_{\mathrm{eff}}}+\left\langle\frac{k\alpha_1}{n}\right\rangle\rfloor}\frac{\Gamma\left(ds_1+\left\langle\frac{k\alpha_2}{n}\right\rangle +\frac{d}{n_{\mathrm{eff}}}\right)}{\Gamma\left(ds_1-\left\langle\frac{k\alpha_1}{n}\right\rangle +1\right)}
\end{align}
For $k$ and $d$ satisfying (\ref{degtwist}), the disc function is then:
\begin{equation}
\label{discfunction}
D_k\left(d, a\right)=\frac{1}{n_{in}} \left(\frac{\hbar}{d}\right)^{\mathrm{age}(\varepsilon^k)} \frac{1}{\lfloor \frac{d}{n_{\mathrm{eff}}}\rfloor!}\frac{\Gamma\left(ds_1+\left\langle\frac{k\alpha_2}{n}\right\rangle +\frac{d}{n_{\mathrm{eff}}}\right)}{\Gamma\left(ds_1-\left\langle\frac{k\alpha_1}{n}\right\rangle +1\right)}
\end{equation}

 
\subsubsection{The localization formula for open invariants} 
\label{sec:locformula}
We combine all ingredients and write down our localization formula/definition for open Gromov-Witten invariants. 
\begin{defn}
For an invariant  for a genus $g$ bordered Riemann Surface with $r$ (labeled) boundary components with winding $d_1,\ldots,d_r$ and $m_i$ insertions of the inertia class $\mathbf{1}_{\frac{i}{n}}$ (and at least two total insertions), we have:
\begin{equation}
\label{locformula}
\begin{array}{l}
{\langle \mathbf{1}^{m_0} \mathbf{1}_{\frac{1}{n}}^{m_1}\ldots \mathbf{1}_{\frac{n-1}{n}}^{m_{n-1}}\rangle_g^{d_1,\ldots,d_r}=}  \\ \\
{\left(\frac{n_{eff}}{\hbar}\right)^r\sum_{k_j} \prod_{j=1}^r \left(\delta^j_0\delta^j_1\delta^j_2 D_{k_j}\left(d, a\right)\right)
\int_{\mathcal{M}}\frac{e^{\mathrm{eq}}(\IE^\vee_{\alpha_0}(s_0)\oplus\IE^\vee_{\alpha_1}(s_1)\oplus\IE^\vee_{\alpha_2}(s_2))}{\prod_{j=1}^r \left({\frac{\hbar}{d_j}}-\psi_j\right)}},
\end{array} 
\end{equation}
where 
$$
\mathcal{M}=\overline{\mathcal{M}}_{g,\sum m_j+r}(B\IZ_n,0;\mathbf{1}^{m_0} \mathbf{1}_{\frac{1}{n}}^{m_1}\ldots \mathbf{1}_{\frac{n-1}{n}}^{m_{n-1}},\mathbf{1}_{\frac{n-k_1}{n}},\ldots, \mathbf{1}_{\frac{n-k_r}{n}})
$$
and the sum is over all $0\leq k_j<n$ that satisfy (\ref{degtwist}) and such
that $\sum k_j= \sum im_i$ (mod $n$). 
\end{defn}
\begin{defn}
\label{def:aoppot}
Let $\lambda$, $\mathbf{w}=\l\{w_m\r\}_{m\in\bbN}$, $\mathbf{\tau}=\l\{\tau_i\r\}_{i=1}^n$ be formal
parameters. The {\it open orbifold Gromov-Witten potential} $F^{(\mathfrak{X},
  \LL)}(\lambda, \mathbf{w},\mathbf{\tau})$ and the {\it genus $g$, $h$-holes
  open orbifold Gromov-Witten potentials} $ F_{g,h}^{(\mathfrak{X}, \LL)}(w_1, \dots,
w_h,\mathbf{\tau})$ of
$(\mathfrak{X}
, L)$
are defined as the formal power series
\bea
\label{eq:aopenpot}
F^{(\mathfrak{X}, \LL)}(\lambda, \mathbf{w},\mathbf{\tau}) &:=&
\sum_{g,h=0}^\infty\lambda^{2g-2+h}\sum_{\stackrel{d_1, \dots d_h}{m_0, \dots,
m_{n-1}}}\prod_{j=1}^h \frac{w_j^{d_j}}{d_j!}\prod_{k=1}^n
\frac{\tau_k^{m_k}}{m_k!} \nn \\ & & 
    {\langle \mathbf{1}^{m_0} \mathbf{1}_{\frac{1}{n}}^{m_1}\ldots
      \mathbf{1}_{\frac{n-1}{n}}^{m_{n-1}}\rangle_g^{d_1,\ldots,d_r}} \nn \\
&=:&
\sum_{g,h=0}^\infty\lambda^{2g-2+h} F_{g,h}^{(\mathfrak{X}, \LL)}(w_1, \dots, w_h,\mathbf{\tau})
\eea
\end{defn}
We refer to the potentials $F_{g,h}^{(\mathfrak{X}, \LL)}(w_1, \dots,
w_h,\mathbf{\tau})$ in terms of the topology of the source curve; in
particular, $F_{0,1}^{(\mathfrak{X}, \LL)}(w,\mathbf{\tau})$ and $F_{0,2}^{(\mathfrak{X}, \LL)}(w_1, w_2,\mathbf{\tau})$ will often be respectively called the {\bf disc potential} and the {\bf
  annulus potential} in the following.

\subsubsection{Disc Invariants and Givental's $J$-function}
An immediate consequence of formula (\ref{locformula}) is that a generating function for disc invariants can be obtained by appropriately turning our disc function into an orbifold cohomology valued function, and pairing it with Givental's $J$-function. 
This is the first step of a general philosophy, that should allow to recover a generating function for all open Gromov-Witten invariants for $[\IC_3/\IZ_n]$ in terms of the (full descendant) Gromov-Witten potential for the closed theory. We are investigating this together with Hsian-Hua Tseng. 

Givental's $J$-function is a generating function that encodes all Gromov-Witten invariants with at most one descendant insertion. We consider the ``small" $J$-function, where we  set the age zero and age two insertion variables equal to $0$. We denote by $\mathbf{1}_\alpha$ the fundamental classes of inertia strata of age one, $\tau_\alpha$ the corresponding dual coordinate. By  $\mathbf{1}_\beta$ we denote an arbitrary inertia stratum.

\bea
J(s_i;\tau_{\alpha};z) &=& 
z\mathbf{1}+\tau_{\alpha}\mathbf{1}_{\alpha}  \sum_{m}\frac{\tau_{\alpha}^m
  (\mathbf{1}_{\beta})^\vee}{m!}\nn \\ &\times
&\int_{\smbzen{m+1}{\mathbf{1}_{\alpha}^{m},
    \mathbf{1}_{\beta}}}\hspace{-1.5cm}
\frac{e(\bE_1^\vee\otimes\mathcal{O}(s_0)\oplus\bE_1^\vee\otimes\mathcal{O}(s_1)\oplus
  \bE_1^\vee\otimes\mathcal{O}(s_2))}{(z - \psi)} \nn \\
\eea

Note that inside the summation formula we insert cohomology classes that are
dual to $\mathbf{1}_{\beta}$ with respect to the orbifold Poincar\'e pairing. \\

We package disc functions into a cohomology valued generating function:
\beq
\mathcal{D}(d,a):= \sum_{k=0}^{n-1} D_k(d,a) (\mathbf{1}_{\frac{k}{n}})^\vee.
\eeq

Then the degree $d$ disc potential for $[\IC^3/\IZ_n]$ is obtained by specializing the variable $z$ to $\hbar/d$ and pairing with the disc function:
$$
\mathcal{F}_0^{disc}(x,y,a):=\sum_{n,d}\langle\mathbf{1}_{\alpha}^m;a\rangle^{d}_{0}\frac{x^m}{m!}\frac{y^d}{d!}=
$$

$$
\sum_d\left[\frac{1}{\hbar}J\left(s_1,s_2,s_3;x;\frac{\hbar}{d}\right)\mathcal{D}(d,a) \right]\frac{y^d}{d!}
$$ 

\begin{rem}
Note that the $J$-function packaging takes care of the unstable terms as well:
\begin{description}
	\item[no insertions] these terms are obtained from the multiplication of the term $\frac{1}{\hbar}z\mathbf{1}$ with $D_0(d,a)(\mathbf{1})^\vee$.
\item[one $\mathbf{1}_{\alpha}$ insertion] likewise these terms are obtained as the products
$$
\frac{1}{\hbar}\mathbf{1}_{\alpha} D_1(d,a)(\mathbf{1}_{\alpha})^\vee
$$
\end{description}
\end{rem}

\section{The $B$-model side}
\label{sec:Bmodel}
\subsection{Toric mirror symmetry and spectral curves}

We review the main concepts which lead to the computation of $B$-model generating functions. We first review the mirror symmetry construction of \cite{Hori:2000ck,
Hori:2000kt} of $B$-model mirrors of toric Calabi-Yau threefolds, thereby introducing the notion of mirror spectral curves, as well as its extension to the open string sector \cite{Aganagic:2000gs, Lerche:2001cw, Aganagic:2001nx}. Finally, we review the formalism of \cite{Eynard:2007kz, Bouchard:2007ys} for the computation of open string potentials from the spectral curve, as well as their transformation properties when crossing a wall in the extended K\"ahler moduli space \cite{Aganagic:2006wq, Bouchard:2008gu, Brini:2008rh}.

\subsubsection{A review of closed mirror symmetry}
Let $X$ be a Calabi-Yau threefold. Let $\{\Phi_i\}_{i=1}^{b_2(X)}$  be a basis
of $H_2(X,\bbZ)$ given by
fundamental classes of compact holomorphic curves in $X$, and denote $\{\Phi^i\}_{i=1}^{b_2(X)}$ 
their duals in co-homology. For $\mathbf{t} \in H^{1,1}(X) \simeq H^2(X, \bbC)$,
write $\mathbf{t}= \sum_i t_i \Phi^i$.  
Closed mirror symmetry for Calabi-Yau threefolds (see \cite{MR1677117} for a comprehensive review) turns the computation of the genus zero Gromov-Witten potential of  $X$ 
%
\beq
F_0^X(t)=\frac{1}{3!} (\mathbf{t}, \mathbf{t} \cup \mathbf{t}) + \sum_{\substack{\beta\in H_2(X, \bbZ) \\ \beta \neq 0}}e^{-\mathbf{t \cdot \beta}} N_{0,\beta}
\eeq
where
\beq
(a, b) = \int_X a \cup b \quad \hbox{for } a, b \in H^\bullet(X), \quad  N_{0,\beta}=\int_{[\cM_{0,0}(X,\beta)]^{\rm vir}}1, 
\eeq
%
into the 
computation of periods of the holomorphic
$(3,0)$ form $\Omega$ of a ``mirror'' flat family of Calabi--Yau threefolds
$\hat X\to B$, where $B$ is a complex algebraic orbifold with
$\mathrm{dim}_{\bbC}B=h^{2,1}(\hat X)$:
\beq
t_i(\{a_j\}) =\oint_{A_i} \Omega \qquad \frac{\de \cF_0^{\hX}}{\de t_i}(\{a_j\})=\oint_{B_i} \Omega.
\label{eq:periodsomega}
\eeq
In \eqref{eq:periodsomega}, $\{a_j\}_{j=1}^{h^{2,1}(\hX)}$ are local
co-ordinates on the base $B$, while $\{A_i, B_i\}_{i=1}^n$ are a basis of homology
$3$-cycles $A_i, B_i \in H_3(\hX, \bbZ)$ such that the intersection pairing
has the canonical Darboux form $(A_i,B_j)=\delta_{ij}$, $(A_i, A_j)=(B_i,
B_j)=0$, and canonically fixed by the asymptotic properties of the periods
around a maximally unipotent  monodromy point. The statement of mirror symmetry is then
\beq
 F_0^X(t_i)= \cF_0^{\hX}(a_j(t_i))
\eeq
 \\

By \eqref{eq:periodsomega}, Gromov-Witten invariants of $X$ can be recovered
by explicit knowledge of the periods of the holomorphic $(3,0)$ form of $\hX$,
and therefore of the mirror manifold $\hX$ itself. \\

\begin{figure}[t]
\vspace{1cm}
\begin{minipage}[t]{0.48\linewidth}
\centering
\vspace{0pt} 
\includegraphics[scale=0.3]{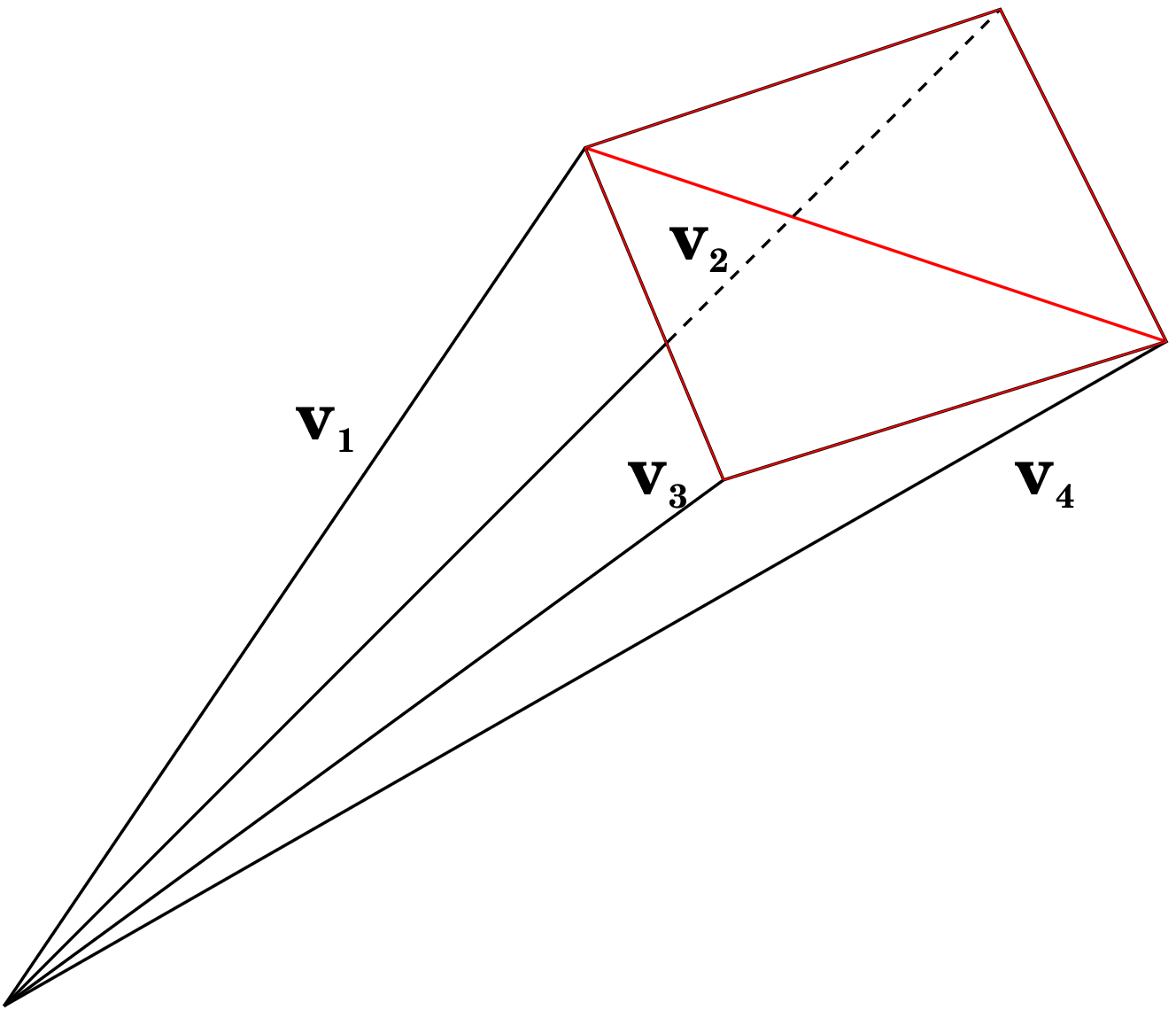}
\vspace{0pt} 
\caption{The toric fan of $X=\cO_{\bbP^1}(-1)\oplus\cO_{\bbP^1}(-1)$.}
\label{fig:fan}
\end{minipage} 
\begin{minipage}[t]{0.48\linewidth}
\centering
\vspace{0pt} 
\includegraphics[scale=0.6]{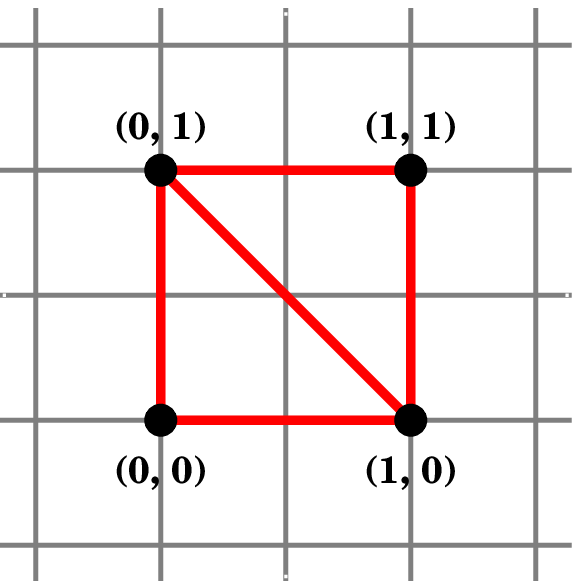}
\vspace{0pt} 
\caption{The Newton polytope $\Sigma_X$ of $X=\cO_{\bbP^1}(-1)\oplus\cO_{\bbP^1}(-1)$.}
\label{fig:sigmax}
\end{minipage} 
\end{figure}

In the case in which $X$ is
toric, it is natural to expect that the pair $(\hX, \Omega)$ could be
constructed entirely from the toric data of $X$. In the physics literature
\cite{Hori:2000kt, Hori:2000ck}, arguments of two-dimensional quantum field
theory suggest an
explicit construction of $(\hat X, \Omega)$, which we briefly review. Since $K_X \simeq \cO_X$, the tip of the $1$-dimensional cones of the fan of $X$ all lie on an affine hyperplane  $H\subset\bbC^3$ \cite{MR1234037}; the intersection of the fan with $H$ yields a finite order subset of $\bbZ^2$ (see Fig. \ref{fig:fan}-\ref{fig:sigmax}). Let $\Sigma_X$ denote the convex hull of such a set of points.
\begin{defn}[Hori-Iqbal-Vafa mirror,  \cite{Hori:2000kt, Hori:2000ck}]
 The $B$-model target space $\hX$ mirror to a toric $CY$ three-fold $X$ is the family of hypersurfaces in $\mathbb{C}^2(x_1,x_2) \times (\mathbb{C}^*)^2(U,V)$
\beq
x_1 x_2 = P_X(U,V)
\label{eq:mirror3fold}
\eeq
where $P_X(U,V)$ is the Newton polynomial associated to the polytope $\Sigma_X$
\beq
\label{eq:newtpol} 
P_X(U,V)=\sum_{p\in \Sigma_X} a_p U^{\rm{pr}_1(p)}V^{\rm{pr}_2(p)}
\eeq
and we have denoted by ${\rm{pr}_i}:\Sigma_X \to \bbZ$, $i=1,2$ the canonical projections to the co-ordinate axes of $\bbC^2 \supset \bbZ^2$.
\label{def:horivafa}
\end{defn}

\begin{example}
Let $X=\cO_{\bbP^1}(-1)\oplus\cO_{\bbP^1}(-1)$. The rays of its fan can be
taken to be
\beq
v_1=\left(\bary{c}1 \\ 0 \\ 1 \eary\right), \quad v_2=\left(\bary{c}0 \\ 1 \\ 1 \eary\right), \quad v_3=\left(\bary{c}0 \\ 0 \\ 1 \eary\right), \quad v_4=\left(\bary{c} 1 \\ 1 \\ 1 \eary\right), \quad
\eeq
The affine hyperplane $H$ in this case is the subspace $z=1$ of $\bbC^3(x,y,z)$. The polytope $\Sigma_X$ is depicted in Fig. \ref{fig:sigmax}; the Newton polynomial $P_X(U,V)$ in this case is
$$P_X(U,V)=a_1 + a_2 U + a_3 V + a_4 U V $$
\end{example}

\noindent Eq. \eqref{eq:mirror3fold} suggests that the non-trivial aspects of the complex geometry of $\hX$ be entirely encoded in the affine curve $P_X(U,V)=0$. This is true in particular for period integrals.
\begin{prop}[\cite{MR2282969, Forbes:2005xt}]
Periods of the holomorphic 3-form reduce to periods
\beq
\oint_\gamma d\lambda_X
\label{eq:periodscl}
\eeq
of the 1-differential
\beq
d\lambda_X = \log{V} \hbox{ d} \log{U}
\label{eq:horivafadiff}
\eeq
over 1-cycles $\gamma$ of the {\bf mirror curve} $H_X$ given by the zero locus $P_X(U,V)=0$. Its projectivization $\Gamma_X=\overline{H_X}$ is a smooth projective curve of genus $\mathfrak{g}$, where $\mathfrak{g}$ is equal to the number of internal points of $\Sigma_X$.
\end{prop}
\begin{rem}
\label{rem:noncpt}
The periods of $\Omega$ for compact Calabi-Yau manifolds are usually computed by
solving the associated Picard-Fuchs system.
However when $X$ is toric, and therefore non-compact, the evaluation of $\Omega$ on $H_3(\hX, \bbZ)\simeq H_1(\Gamma_X, \bbZ)$ fails to reproduce a complete set of solutions of the Picard-Fuchs equations \cite{MR1166813, MR1677117}. If $\Delta$ is a choice of a principal branch for the logarithm on $\Gamma_X$, i.e., a disconnected union of real segments on $\Gamma_X$ such that $\log{U}$, $\log{V}$ are single-valued meromorphic functions on $\Gamma_X \setminus \Delta$, the missing period integrals can be recovered by considering periods of $d\lambda_X$ along non-compact cycles $\gamma$ in $\Gamma_X \setminus \Delta$ \cite{Brini:2008rh, Forbes:2005xt}. When we need to stress that we refer to the set of non-compact periods with logarithmic singularities (i.e. periods over three-cycles which are mirror of non-compact divisors of $X$), we denote them with a tilde $\{\tilde{t}_i\}_{i=1}^{b_2(X)-\mathfrak{g}}$.
\end{rem}

\subsubsection{Open string mirror symmetry}
We have seen that the ordinary statements of mirror symmetry simplify, in the toric case, into computations of periods of a $1$-differential on a Riemann surface. This situation generalizes to the open string setting. \\

Open string mirror symmetry deals with a $B$-model construction of the open
Gromov-Witten potentials $F_{g,h}^{(X,L)}$ (as in definition \ref{def:aoppot})
of a pair $(X,L)$, with $L\subset X$ a Lagrangian submanifold,
%
in terms of a ``mirror'' pair $(\hX, \widehat{L})$, where $\widehat{L}$ is a
holomorphic submanifold of $\hX$.
 As a natural extension of  the closed mirror symmetry lore, genus zero open
mirror symmetry intends to recover 
genus zero {\it open} Gromov-Witten invariants $N_{0, \beta}^{d_1, \dots,
  d_h}$ and the corresponding potentials
\beq
F_{0,h}^{(X, L)}(t_i, x_i) = \sum_{\beta \in H_2(X, \bbZ)} \sum_{d_1=1, \dots, d_h=1}^\infty N_{0, \beta}^{d_1, \dots,
  d_h}
e^{-\mathbf{t} \cdot \beta}
\prod_{i=1}^h\frac{ x_i^{d_i}}{d_i! }
\label{eq:Xopenpot}
\eeq
 from the study of complex variations of the {\it pair} ($\hX$, $\widehat{L}$), thus leading to period computations in {\it relative} co-homology \cite{Lerche:2002yw, Forbes:2003ki, MR2481273}. In particular, the disc potential  (in physical terms, the domain wall tension) should be computed as a co-chain integral \cite{Witten:1997ep}
\beq
\cF_{0,1}^{(\hX,\widehat{L})} = \int_H \Omega 
\label{eq:cochain}
\eeq
where $\de H = B^+-B^-$ and $[B^+]=[B^-]=[\widehat{L}]$.\footnote{Notice that \eqref{eq:cochain} implies the choice of representatives in the homology class of $\widehat{L}$. This causes an ambiguity in the leading term of the open string moduli expansion, entirely analogous to the quadratic ambiguity of the ordinary, closed genus zero Gromov-Witten potential.} \\

The toric case presents a number of simplifications in the open setting too. A
distinguished class of special Lagrangian $A$-branes  with topology
$\bbR^2\times S^1$  were constructed by Aganagic and Vafa in
\cite{Aganagic:2000gs}: in an affine patch, these are  the Lagrangians
constructed in Sect. \ref{sec:Amodel} for $n=1$. It was proposed in
\cite{Aganagic:2000gs} that their mirror $B$-branes should be cut by the equations
\beq
x_1 = 0 = P(U,V) \qquad \hbox{or} \qquad x_2 = 0 = P(U,V)
\eeq
The ambiguity in the choice of $x_1$ or $x_2$  results \cite{Aganagic:2000gs}  in an overall sign ambiguity of the open string amplitudes. For this kind of branes, dimensional reduction of the holomorphic Chern-Simons action on the brane shows that the computation of disc invariants reduces to the computation of a sort of ``Abel-Jacobi'' map on the mirror curve.
\begin{defn}
Let $X$ be a toric Calabi-Yau threefold, and $L\subset X$ be an Aganagic-Vafa
Lagrangian $A$-brane. Then the {\rm $B$-model disc potential} of ($\hX$, $\widehat{L}$) is given by
\beq 
\cF_{0,1}^{(\hX,\widehat{L})}(p) = \int_{p_*}^{p} d\lambda_X
\label{eq:bdisc}
\eeq
where $p_*$, $p \in \Gamma_X$, with $p_\ast$ fixed, and $d\lambda_X$ is as in \eqref{eq:horivafadiff}.
\end{defn}
%
\begin{rem}
 Both the ``closed'' \eqref{eq:periodscl} and the ``open'' \eqref{eq:bdisc}
 periods are defined in terms of a contour integral of the one-form
 $d\lambda_X$, which is specified by the toric data. The latter in itself is
 however only defined up to an action of $\cG \simeq GL(3,\bbZ)$, i.e. changes
 of basis for the three dimensional lattice where the fan of $X$ lives; in particular, a subgroup $\HH \simeq GL(2,\bbZ)$ acts effectively on the hyperplane $H$ where the tip of the 1-dimensional cones lie.  By \eqref{eq:newtpol}, this induces a $GL(2,\bbZ)$ transformation on the $B$-model variables $V$ and $U$
\bea
\label{eq:gl2z}
\left(\bary{c}V \\ U
\eary\right) & \rightarrow &
\left(\bary{c}\tilde V \\ \tilde U
\eary\right)=\left(\bary{c}V^a U^b \\ V^c U^d
\eary\right), \quad ad-bc=\pm 1 \\
\eea
and, accordingly, on the 1-differential $d\lambda_X=\log{V} d\log{U}$. Remarkably
\cite{Eynard:2007kz}, the prepotential $\cF_0^{\hX}$ computed via
\eqref{eq:periodsomega} is invariant under the transformation
\eqref{eq:gl2z}. This is however not the case for the disc potential
\eqref{eq:bdisc}: that is, the disc potential is not an invariant of the pair ($\hX$,$\widehat{L}$), but it rather comes with an integer ambiguity. Its meaning was elucidated in \cite{Aganagic:2001nx} (see also \cite{Bouchard:2007ys} for a very clear exposition). First of all recall that $GL(2,\bbZ)$ has three generators:

\beq
P= \left(\bary{cc}1 & 0 \\ 0 & -1
\eary\right), \quad 
S= \left(\bary{cc}0 & -1 \\ 1 & 0
\eary\right), \quad
T= \left(\bary{cc}1 & 1 \\ 0 & 1
\eary\right)
\eeq
The $T$ transformation generates a free abelian subgroup of $SL(2,\bbZ)$ which leaves invariant the $U$-direction of $\bbC^2(V,U)$. Then:
\ben
\item Fixing a $U$-direction, for example by acting by a combination  involving $S$ and $T$, amounts to specifying the Aganagic-Vafa SLag $L$. The reader can find the details, based on the description of smooth toric Calabi-Yau threefolds as $\bbT^2 \times \bbR$ fibrations, in \cite{Bouchard:2007ys}. 
\item After fixing the $U$-direction, there's a leftover $\bbZ \rtimes
  \bbZ_2$-ambiguity given by the $T$ and $P$ transformations. The
  $P$-ambiguity results by \eqref{eq:horivafadiff} in a sign ambiguity in the
  definition of the disc function, presumably related to orientation problems
  in the construction of the moduli space of stable maps with Lagrangian
  boundary conditions \cite{Solomon:2006dx}. The $T$-ambiguity, called the
  {\it framing} of the brane $L$, is an intrinsic ambiguity in the computation
  of the disc function, and it was related in  \cite{Aganagic:2001nx} to an
  analogous ambiguity \cite{Witten:1988hf} in the conjectural dual description
  of the $A$-model on $(X,L)$ via Chern-Simons theory and related knot
  invariants \cite{Ooguri:1999bv}. 
\een 
\end{rem}
\begin{conj}[Mirror symmetry for disc invariants]
\beq
F_{0,1}^{(X, L)}(t_i, x_j) = 
\cF_{0,1}^{(\hX, \widehat{L})}(z_i(t_i), p(t_i,x_j))
\eeq
\end{conj}

As in the ordinary closed string case, physical arguments
  related to the BPS interpretation of open string amplitudes suggest that
  the conjectural relationship of $\cF_{0,1}^{\hX, \widehat{L}}$ with a
  Gromov--Witten disc potential should hold true \cite{Aganagic:2001nx,
    Lerche:2001cw} only up to a change of variables relating the $B$-model
  open modulus $p$ in \eqref{eq:bdisc}, i.e. a point on the mirror curve
  $\Gamma_X$, with a suitably defined $A$-model open co-ordinate $x$. 
Mathematically, this is achieved by writing a Picard-Fuchs system extended to relative co-homology \cite{Forbes:2003ki}: the additional solutions provide the so-called {\it open mirror map}. When $X$ is smooth, i.e. at ``large radius'',  $p$ in \eqref{eq:bdisc} and $x$ in \eqref{eq:Xopenpot} are related as
\beq
\label{eq:openflat}
x = p \prod_i \left(\frac{z_i}{q_i}\right)^{r_i}
\eeq
where $r_i \in \bbQ$, $z_i$ are $B$-model closed moduli, and $q_i$ are
exponentiated, closed flat co-ordinates $q_i=e^{t_i}$; the rational
numbers $r_i$ are determined by the solutions of the extended Picard-Fuchs
system. In other words, Eq. \eqref{eq:openflat} means that the open string $A$--model modulus is related to the $B$--model one by a correction involving \textit{closed} moduli only. Equation (\ref{eq:openflatA}) describes how \eqref{eq:openflat} is modified in the orbifold setting.

\subsection{The remodeled $B$-model and open orbifold invariants}
\subsubsection{The Eynard-Orantin recursion}
We have seen how the  $B$-model prepotential and disc function are completely determined in terms of the mirror geometry, i.e., the mirror curve $\Gamma_X$ together with its graph in $(\bbC^*)^2(U,V)$. Recently, an influential proposal was put forward by  Bouchard, Klemm, Mari\~no and Pasquetti \cite{Marino:2006hs, Bouchard:2007ys}, which gives a complete and unambiguous prescription for the computation of generating functions for genus $g$, $h$-holed open Gromov-Witten invariants via residue calculus on $\Gamma_X$. Their conjecture was based on an application of the Eynard-Orantin recursive formalism \cite{Eynard:2007kz} to the case of mirrors of toric Calabi-Yau threefolds. \\
 To give the precise statement of the conjecture, we start with the following

\begin{defn}
\label{defn:spectral}
A spectral curve $\cS$ is a 5-tuple $(\Gamma, C, \Delta, u,v)$ where
\ben
\item $\Gamma$ is a family of genus $\mathfrak{g}$  complex projective curves, 
\item $C=\{C_1, \dots, C_m\}$, for $m\in \bbN$, is a collection of holomorphic sections of $\Gamma$,
\item $\Delta=\{\Delta_1,  \dots, \Delta_{[\frac{m}{2}]}\}$ is a smooth real
  family of arcs $\Delta=\{\Delta_i = (C_{2i-1},C_{2i})\}_{i=1}^{[\frac{m}{2}]}$, 
\item $u,v:\Gamma\to \bbC$ are marked analytic functions on $\Gamma$, meromorphic on $\Gamma \setminus \Delta$ and with at most logarithmic polydromies on 
$\Delta$.\een

If $du$ and $dv$ never vanish simultaneously,
the spectral curve is called {\it regular}. 
\end{defn}

Mirror symmetry for toric  Calabi-Yau threefolds provides us with an example
of a spectral curve. In this case $\cS=(\Gamma_X, C, \Delta, \log{U},
\log{V})$, where $C=\{ p \in \Gamma_X | V(p) U(p)=0\} \cup \{ p \in \Gamma_X |
1/(V(p) U(p))=0\}$, $\Delta$ is a choice of principal branch for $\log V$,
$\log U$, and we have denoted with the same symbol $U$, $V$ the unique meromorphic lift of $U$, $V$ to $\Gamma_X$. \\

Suppose now that $\cS$ is a regular spectral curve, and let $\{q_i\}$ denote the ramification points of the $v$ projection to $\bbC$. Near  $q_i$ there are two points $q, \bar q \in \Gamma$ with the same projection $U(q)=U(\bar q)$. Picking a polarization $\HH \in Sp(2\mathfrak{g}, \bbZ)$ of $\Gamma$, that is a symplectic basis of $H_1(\Gamma, \bbZ)$, the {\it Bergmann kernel} is defined as the unique meromorphic differential on $\Gamma\times\Gamma$ with a double pole at $p=q$ with no residue and no other pole,
and normalized so that for every $p\in \Gamma $
\beq
\label{eq:bergmannnormalization}
\oint_{p\times A_I} B(p,q) = 0,
\eeq
It is useful to introduce also the $1$--form
\beq
\label{deqp}
d E_{q} (p)= {1 \over 2} \int_{q}^{\bar q} B(p,\xi),
\eeq
which is defined locally near a ramification point $q_i$. Notice that $B(p,q)$
depends only on $(\Gamma, \HH)$ and on no additional data. \\ 

Out of $\cS$, Eynard and Orantin  \cite{Eynard:2007kz}  define recursively an
infinite sequence of {\it correlators} $W_h^{(g)} (p_1, \ldots, p_h)$
from the spectral curve as follows:

\begin{defn}[Eynard--Orantin recursion]
For all  $g,h \in \IZ^+$, $h \geq 1$, a meromorphic differential $W_h^{(g)} (p_1, \ldots, p_h)\in \mathrm{Sym}^h \Omega^{(1,0)}(\Gamma_X)$  is defined from the following recursion

\bea
\label{eq:BKMPrecursion1}
W_1^{(0)}(p) &=& 0\\
W_2^{(0)}(p,q) &=& B(p,q) \\
W_{h+1}^{(g)}(p, p_1 \ldots, p_h) &=& \sum_{q_i}  \underset{q=q_i}{\rm Res~} {d E_{q}(p) \over \Phi(q) - \Phi(\bar q)} \Big ( W^{(g-1)}_{h+2} (q, \bar q, p_1, \ldots, p_{h} )\nn \\
&+& \sum_{l=0}^g \sum_{J\subset H} W^{(g-l)}_{|J|+1}(q, p_J) W^{(l)}_{|H|-|J| +1} (\bar q, p_{H\backslash J}) \Big)
\label{eq:BKMPrecursion3}
\eea

Here we wrote $\Phi(u) := v(u) du$. Moreover we denoted $H={1, \cdots, h}$, and given any subset $J=\{i_1, \cdots,
i_j\}\subset H$ we defined $p_J=\{p_{i_1}, \cdots, p_{i_j}\}$. 
\end{defn}

The entire set of correlators
is constructed out of the spectral curve by residue calculus on
$\Gamma$. The conjecture of \cite{Marino:2006hs, Bouchard:2007ys}
is that, when $\cS$ is the mirror spectral curve of a toric Calabi-Yau
threefold $X$, such quantities compute precisely the open 
Gromov--Witten potentials of $(X,L)$, for any genus $g$ and number of holes $h$. 

\begin{conj}[BKMP, \cite{Marino:2006hs, Bouchard:2007ys}]
\label{conj:BKMP}
Let $\cS=(\Gamma_X, C, \Delta, \log{U},
\log{V})$ be the mirror spectral curve of a toric $CY$ $3$--fold $X$, and let
$A_i$ in (\ref{eq:bergmannnormalization}) correspond to homology $1$-cycles in
$\Gamma_X$ such that the periods of the Hori--Vafa differential
have logarithmic singularities at the large complex structure point. 
Let $\cS_f$ be the one--integer parameter family of spectral curves obtained
by sending $U\to UV^f$, $V\to V$ for $f\in \bbZ$. Then the integrated
correlation functions
$\mathcal{F}^{(\widehat{X},\widehat{L})}_{g,h} = \int
W_k^{(g)}(p_1,\ldots,p_k) \frac{dp_1}{p_1} \dots \frac{dp_h}{p_h}$, for $2g+h>1$, are equal to the A-model framed open Gromov--Witten potential of $(X,L)$ where $L$ is the mirror brane to $\Gamma_X\subset \hX$, after plugging in the closed and open mirror maps. 
\end{conj}

\subsubsection{Open string amplitudes and wall-crossings}
\label{sec:BKMP2}
The residue computation of
Eqns. \eqref{eq:BKMPrecursion1}-\eqref{eq:BKMPrecursion3} gives, in principle,
the correlators $W_h^{(g)} (p_1, \ldots, p_h)$ as closed functions of the open
moduli $p_1,\dots p_h$ as well as of the complex moduli of the Hori-Iqbal-Vafa
curve \eqref{eq:newtpol}. A remarkable property of $W_h^{(g)} (p_1, \ldots,
p_h)$ is that they are {\it almost-modular forms} \cite{Eynard:2007kz,
  Bouchard:2008gu} of $\Gamma_X$, as we now review. \\

The mirror Calabi-Yau $\hX$ of $X$ has a complex structure moduli space
$\cM_\hX$, which by the Bogomolov-Tian-Todorov theorem is a smooth complex
manifold of complex dimension $b_2(X)$. $\cM_\hX$ admits a natural toric compactification to a toric orbifold
$\overline{\cM_\hX}$, whose fan is given by the 
secondary fan of $X$ \cite{MR1677117}; the Gauss-Manin connection on $\cM_\hX$ lifts to an (in general meromorphic) connection on $\overline{\cM_\hX}$, whose monodromies around each boundary point of $\cM_\hX$ generate the {\it monodromy group} $\cG$ of $\hX$. The latter \cite{Aganagic:2006wq} turns out to be a finite index subgroup of $Sp(2\mathfrak{g},\bbZ)$, where $\mathfrak{g}$ is the genus of $\Gamma_X$. We have the following
\begin{thm}[\cite{Eynard:2007kz, Bouchard:2008gu}]
$W_h^{(g)} (p_1, \ldots, p_h)$ admits the following expansion
\beq
W_h^{(g)} (p_1, \ldots, p_h) = \sum_{n=0}^{3g-3+2h} c_{n}(\tau, \{\tilde{t_i}\},\{p_j \}) E_2^{n}(\tau)
\label{eq:WhgE2}
\eeq
where $\tau$ is the period matrix of $\Gamma_X$,
$\{\tilde{t_i}\}_{i=1}^{b_2(X)-\mathfrak{g}}$ are the periods of $v(p)du(p)$
over cycles mirror of non-compact divisors of $X$ (see Remark
\ref{rem:noncpt}), $c_{n}$ is a holomorphic function of $\tau$,
$\{\tilde{t_i}\}$, and $\{p_j \}$, and $E_2$ is the genus-$\mathfrak{g}$ generalization of the second Eisenstein series (see \cite{Aganagic:2006wq}).
\label{thm:modular}
\end{thm}
Theorem \ref{thm:modular} acquires particular relevance in view of the following
\begin{conj}[\cite{Bouchard:2007ys}]
\label{conj:modular}
$W_h^{(g)} (p_1, \ldots, p_h)$ is a weight zero holomorphic almost modular form of $\cG$. More precisely, $c_{n}$ in \eqref{eq:WhgE2} is for all $p_i$ and $\tilde t_i$ a $-2n$ modular form of $\cG$.
\end{conj}
Let us explain more in detail what we mean by almost modularity, focusing for definiteness to the case $\mathfrak{g}=1$ which will be discussed in Sect. \ref{sec:c3z3}-\ref{sec:z4} . Under an $Sp(2,\bbZ)=SL(2,\bbZ)$ transformation
\bea 
\label{eq:msl2z}
M &:=& \left(\bary{cc}A & B \\ C & D  \eary\right) \in SL(2,\bbZ)  \\
\tau & \to & \tilde \tau = (C \tau+ D)^{-1}(A \tau + B)
\label{eq:modtrasf}
\eea
$E_2(\tau)$ transforms as
\beq
E_2(\tilde \tau)=  (C \tau+ D)^2 E_2(\tau)+ d_2(\tau)
\label{eq:e2trasf}
\eeq
where
\beq
d_2(\tau)=\frac{6}{\pi i }C(C \tau + D)
\label{eq:dtau}
\eeq
Hence, it is nearly a weight two modular form of $SL(2,\bbZ)$, but for a shift linear in $\tau$. The almost modularity of $W_h^{(g)} (p_1, \ldots, p_h)$ stems entirely from that of $E_2(\tau)$. Under a modular transformation $\tau \to \tilde \tau$, the expansion 
\eqref{eq:WhgE2} gets transformed to
\bea
W_h^{(g)} &\to & \cW_h^{(g)} \nn \\
\cW_h^{(g)} (p_1, \ldots, p_h) &=& \sum_{n=0}^{3g-3+2h} c_{n}(\tau, \{\tilde{t_i}\},\{p_i \}) (E_2(\tau)+d_2(\tau))^n
\label{eq:WhgE2trasf}
\eea
Eq. \eqref{eq:WhgE2trasf} expresses the variation of the open string generating functions under a change in the choice of polarization of the mirror curve. Recall that in Conjecture \ref{conj:BKMP} a polarization was fixed by requiring the $A$--periods of the Hori--Vafa differential to be large radius flat co-ordinates, i.e. logarithmic solutions of the $PF$ system around the maximally unipotent monodromy point. Changing polarization then corresponds to an $Sp(2\mathfrak{g}, \bbZ)$ transformation to a different basis of solutions of the $GKZ$ system. 
\\

\noindent Almost-modularity has a particular relevance for the mirror symmetry
treatment of the behaviour of Gromov--Witten potentials under variations of
the K\"ahler structure, and in particular under birational
transformations. Let us consider a situation in which two pairs $(X, L)$ and
$(\mathfrak{X}, \LL)$ are given, where $X$ is a smooth toric CY3, $L \subset X$
an Aganagic-Vafa brane, $\mathfrak{X}$ a reduced algebraic orbifold birational
to $X$ and $\LL$ the corresponding lagrangian in $\mathfrak{X}$. Let
$T\simeq\bbC^* \circlearrowright X, {\mathfrak{X}}$ specify torus actions on $X$, ${\mathfrak{X}}$ that  act trivially on the canonical bundle and such that the resolution morphism is $T$-equivariant.
 We  use here $\{t^X_i\}_{i=1}^{b_2(X)}$ and $\{t^{\mathfrak{X}}_i\}_{i=1}^{b^{\mathrm{CR}}_2({\mathfrak{X}})}$ for the quantum parameters of
 $QH_T^\bullet(X)$ and $QH_T^\bullet(\mathfrak{X})$ respectively\footnote{In
   Section~\ref{sec:locformula} we used $\tau_i$ for the variables of orbifold
 quantum co-homology; we prefer to denote them with $t^{\mathfrak{X}}_i$ here to
 avoid confusion with the period matrix $\tau$ of the mirror curve.}, and
 $\{x_j\}$ and $\{ w_j\}$ for their open string expansion
 parameters. In the terminology of  Remark  \ref{rem:noncpt}, we distinguish between
 ``compact'' moduli $\{\hat t^X_{i}, \hat t^{\mathfrak{X}}_{i} \}_{i=1}^{\mathfrak{g}}$
 and ``non-compact'' ones $\{\tilde t^X_{i}, \tilde t^{\mathfrak{X}}_{i}
 \}_{i=1}^{b_2(X)-\mathfrak{g}}$. Mirror symmetry arguments then lead to the following statements:
\ben 
\item  the (compact) flat coordinates and the prepotential \cite{MR2510741, MR2454327} of $X$ and $\mathfrak{X}$ should be related by a linear, $2\mathfrak{g} \times 2\mathfrak{g}$ invertible transformation
\beq
\left(\bary{c} \hat t^X_i \\  \frac{\de \cF_0^X}{\de \hat t^X_i}
\eary\right)=
\left(\bary{cc}
A & B  \\
C & D 
\eary\right)
\left(\bary{c} \hat t^{\mathfrak{X}}_i \\  \frac{\de \cF_0^{\mathfrak{X}}}{\de
  \hat t^{\mathfrak{X}}_i}
\eary\right)
\label{eq:M2}
\eeq
\item the winding parameters  $\{x_j\}$ and $\{w_j\}$ should be related by a rescaling factor, involving exponentiated flat coordinates only
\beq
x_j=  w_j \prod_k q_k^{r_k}
\label{eq:openflatA}
\eeq
\een

We refer to Eq. (51) as the {\it open orbifold mirror map}. 
\begin{rem}
Eq. \eqref{eq:M2} was justified physically in \cite{Aganagic:2006wq} as a necessary transformation to ensure  monodromy invariance of the orbifold partition function. Its ultimate mathematical justification resides in Givental's symplectic vector space formalism \cite{MR2454327, MR2486673}. Eq. \eqref{eq:openflatA} was taken in \cite{Bouchard:2007ys, Brini:2008rh} as a working definition of an open ``flat'' modulus at the orbifold point: in the examples of \cite{Bouchard:2007ys, Brini:2008rh}, this was the minimal choice that could yield an analytic potential at the orbifold point, without fractional powers of the quantum parameters. An {\it a priori} derivation of \eqref{eq:openflatA}, even by physics-based considerations, is to our knowledge still lacking. \\
\end{rem}

The following conjecture states that knowledge of the $\mathfrak{g} \times \mathfrak{g}$ matrices $A$, $B$, $C$ and $D$ suffices to reconstruct the open Gromov-Witten potentials of $\mathfrak{X}$ starting from those of $X$.

\begin{conj}[Mirror symmetry for open orbifold invariants]
\label{conj:BKMP2}
Let $W_h^{(g)}$ denote the open string correlators  of $X$ for $2g+h>1$; when
$g=0$, $h=1$ define $W_1^{(0)}=d\lambda_X$. Let moreover $M$ be the matrix
\beq
\left(\bary{cc}
A & B  \\
C & D 
\eary\right)
\eeq
representing the change of basis from the (normalizable) solutions of the $PF$ system at large radius to those of the $B$--model boundary point associated to $\mathfrak{X}$. 
Define now for $2g+h>1$ the transformed open string correlators $\cW_h^{(g)}$ of
$\mathfrak{X}$ as in \eqref{eq:WhgE2trasf}; when $g=0$, $h=1$, set $\tilde
W_1^{(0)}=W_1^{(0)}$. Then the open orbifold potentials
$F_{g,h}^{(\mathfrak{X}, \mathcal{L})}$ in \eqref{eq:aopenpot} are given by the integrated correlator $\int \cW_k^{(g)}(p_1,\ldots,p_k)\frac{dp_1}{p_1} \dots \frac{dp_h}{p_h}$, after plugging in the orbifold open and closed mirror maps.
\end{conj}
\noindent Conjecture \ref{conj:BKMP2} thus prescribes a three-steps recipe to compute open Gromov--Witten invariants of $\mathfrak{X}$ starting from those of $X$:
\ben
\item when $2g+h>1$, transform the correlators as in \eqref{eq:WhgE2trasf};
\item analytically continue them from the large radius to the relevant boundary point corresponding to $\mathfrak{X}$;
\item expand them in powers of the appropriate local flat co-ordinates.
\een
\begin{rem}
It should be noticed that the two bases of solutions of the $PF$ system need not be related by a simple change of polarization of the mirror curve. This is particularly true for the case of orbifolds \cite{Aganagic:2006wq, Brini:2008rh}. In that case, however, Eqs. \eqref{eq:dtau} and \eqref{eq:WhgE2trasf} still make sense, even though they are no longer the result of the composition of $W_{h}^{(g)}$ with the modular transformation \eqref{eq:modtrasf}. This is why Eq. \eqref{eq:WhgE2trasf} was taken as the definition of the transformed $\cW_h^{(g)}$ in Conjecture \ref{conj:BKMP2}.
\end{rem}
\section{Warming up: $[\bbC^3/\bbZ_3]$}
\label{sec:c3z3}
In this section we specialize the computation to the case of disc invariants
for the orbifold $[\IC^3/\IZ_3]$. We first review the $B$-model predictions by
Bouchard, Klemm, Mari\~no and Pasquetti  in \cite{Bouchard:2007ys,   Bouchard:2008gu}, and then recover them via our formalism. A similar
computation was carried out independently by Hsian-Hua Tseng \cite{ht:pc}.

\subsection{The $B$-model disc potential}

\begin{figure}[h]
\begin{minipage}[t]{0.48\linewidth}
\centering
\vspace{0pt} 
\includegraphics[scale=0.5]{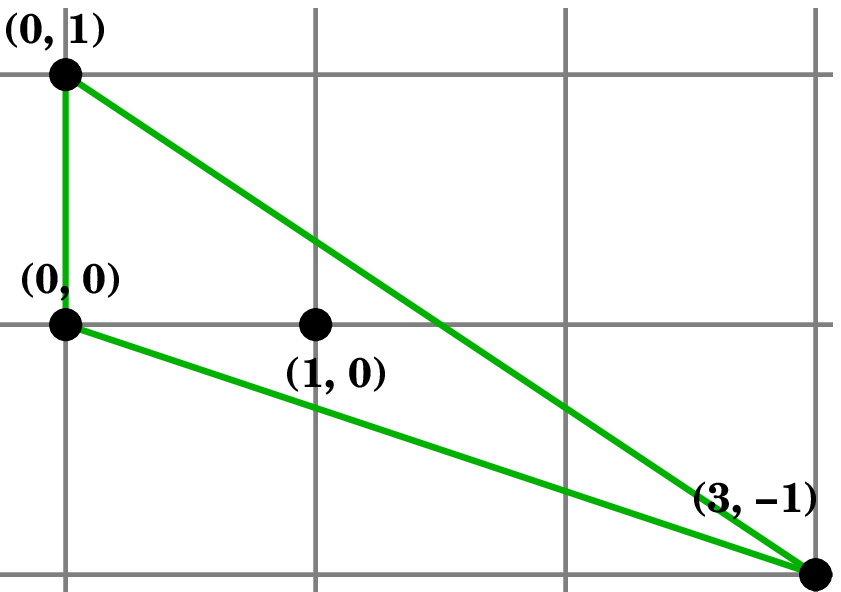}
\vspace{0pt} 
\caption{The Newton polytope of the fan of
  $\mathfrak{X}=[\bbC^3/\bbZ_3]$}.
\label{fig:c3z3}
\vspace{1.5cm}
\end{minipage} 
\begin{minipage}[t]{0.48\linewidth}
\centering
\vspace{0pt} 
\centering
  \includegraphics[scale=0.5]{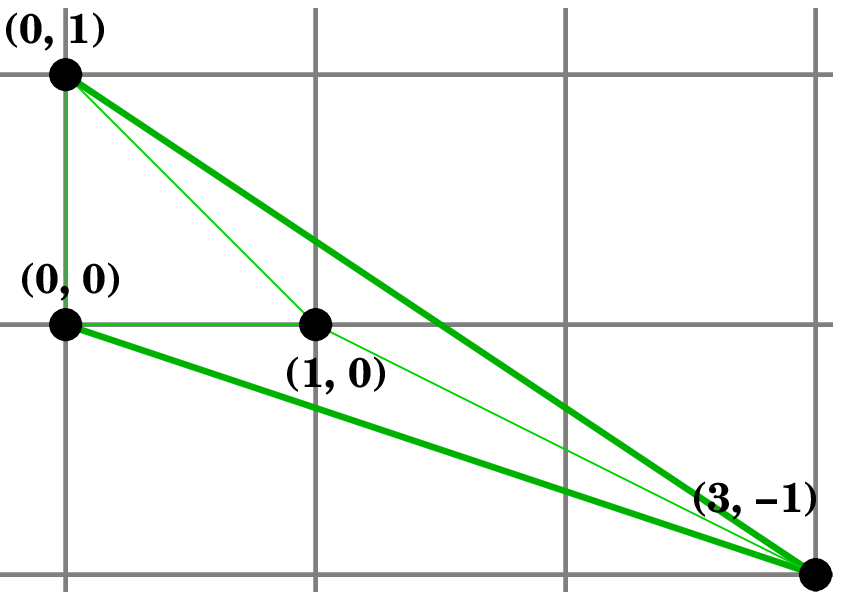}
\caption{The Newton polytope of the fan of $X=K_{\bbP^2}$.}
\label{fig:kp2}
\end{minipage}
\end{figure}

The orbifold $[\IC^3/\IZ_3]$ is obtained by quotienting affine space with
characters $\alpha_0=\alpha_1=\alpha_2=1$; the Newton
polytope associated to its fan is represented in Figure \ref{fig:c3z3}. The crepant resolution is the
canonical line bundle over the projective plane $K_{\bbP^2}$ (Newton
polytope in Figure \ref{fig:kp2}).

According to Definition \ref{def:horivafa}, the mirror curve of
$[\IC^3/\IZ_3]$ is given by
\beq
1 + V + U- \frac{U^3}{3 \psi V}=0
\label{eq:mirrorcurvez3}
\eeq
where $\psi$ is the $B$-model mirror of the $A$-model flat co-ordinate $\tau$,
and we write $QH^2([\IC^3/\IZ_3]) \ni \Phi = \tau_{\frac{1}{3}} \mathbf{1}_{1/3}$ for a
generic age $1$ twisted class. It was argued in \cite{Bouchard:2007ys} that
this choice of representative of the mirror curve corresponds to a brane with
zero framing, located on the outer legs of the $p$-$q$ web diagram of
$K_{\bbP^2}$. Moreover, the authors of \cite{Aganagic:2006wq, Bouchard:2007ys}
found that the mirror map relating $\psi$ and $\tau_{\frac{1}{3}}$ has the form
$$ \tau_{\frac{1}{3}} =  \psi \; {}_3F_2\left(\frac{1}{3}, \frac{1}{3}, \frac{1}{3}; \frac{2}{3},\frac{4}{3}; -\frac{\psi^3}{27}\right)$$
where ${}_3F_2\left(a, b, c; d,e; x\right)$ is the generalized hypergeometric function
$${}_3F_2\left(a, b, c; d,e; x\right) = \sum_{n=0}^\infty \frac{\Gamma(a+n) \Gamma(b+n) \Gamma(c+n) \Gamma(d) \Gamma(e)}{\Gamma(a) \Gamma(b) \Gamma(c) \Gamma(d+n) \Gamma(e+n)}\frac{x^n}{n!} $$
while the open orbifold mirror map is
$$w = \psi U $$
The $B$-model orbifold disc potential at framing zero is then
\beq
\cF^{[\bbC^3/\bbZ_3], (f=0)}_{0,1}(\tau_{\frac{1}{3}},w)=\int^{U(w,\psi(\tau_{\frac{1}{3}}))} d\lambda^{f=0}_{[\bbC^3/\bbZ_3]}(\psi(\tau_{\frac{1}{3}}),U)
\label{eq:Bdiscpotz3}
\eeq
where the Hori-Vafa differential reads, from \eqref{eq:mirrorcurvez3},
$$d\lambda^{f=0}_{[\bbC^3/\bbZ_3]}(\psi,U) = \frac{\log \left(\frac{1}{2}
  \left(\sqrt{(U+1)^2-\frac{4 U^3}{\psi }}+U+1\right)\right)}{U} $$
The inclusion of framing can then be accomplished \cite{Bouchard:2008gu} through the
$T$-transformation 
$$(V,U) \to (V, U V^f) $$

\subsection{The $A$-model disc potential}

 In this case we have only one age $1$ class in orbifold cohomology, namely the class $\mathbf{1}_{\frac{1}{3}}$. Non-equivariant invariants only admit non-trivial insertions of this type. Condition (\ref{degtwist}) and the monodromy condition on the space of maps to $B\IZ_3$ imply degree, twisting and number of insertions are all equal mod $3$ ($d\equiv k\equiv m$ (mod $3$)). Then only one fixed locus contributes to the disc invariant $\langle \mathbf{1}_{\frac{1}{3}}^m\rangle_{0}^{d}$, and formula (\ref{locformula}) reduces to:

\begin{equation}
\label{c3z3disc}
\langle \mathbf{1}_{\frac{1}{3}}^m\rangle_{0}^{d}=
\left\{
\begin{array}{ll}
\frac{1}{d}D_0(d,a) & m=0, d\equiv 0 \ \mbox{(mod $3$)}\\
\frac{1}{\hbar} D_1(d,a) & m=1, d\equiv 1 \ \mbox{(mod $3$)}\\
3D_k(d,a) f_1(m,d,a) &   d\not\equiv 0 \ \mbox{(mod $3$)} \\
3s_0s_1s_2D_k(d,a) f_2(m,d,a) &   d\equiv 0 \ \mbox{(mod $3$)}
\end{array}
\right.
\end{equation}
where
\bea
f_1(m,d,a) &=&
\int_{\overline{\mathcal{M}}_{0,m+1}(B\IZ_3,0;\mathbf{1}_{\frac{1}{3}}^m,\mathbf{1}_{\frac{-m}{3}})}\frac{e^{\mathrm{eq}}(\IE^\vee_{1}(s_0)\oplus\IE^\vee_{1}(s_1)\oplus\IE^\vee_{1}(s_2))}{{\frac{\hbar}{d}}-\psi} \nn
\\ \\
f_2(m,d,a) &=& \int_{\overline{\mathcal{M}}_{0,m+1}(B\IZ_3,0;\mathbf{1}_{\frac{1}{3}}^m,\mathbf{1})}\frac{e^{\mathrm{eq}}(\IE^\vee_{1}(s_0)\oplus\IE^\vee_{1}(s_1)\oplus\IE^\vee_{1}(s_2))}{{\frac{\hbar}{d}}-\psi}
\eea
The torus weights are
$$
(s_0,s_1,s_2)=\left(\frac{1}{3}, -a, a-\frac{1}{3}\right)
$$
and the disc function:
$$
{D}_k(d,a)= \frac{1}{\left\lfloor \frac{d}{3}\right\rfloor!}\left(\frac{\hbar}{d}\right)^{3\langle d/3\rangle}\frac{\Gamma(\frac{d}{3}+\langle \frac{d}{3}\rangle-{da})}{\Gamma(1-\langle \frac{d}{3}\rangle-{da})}
$$

Specializing to the torus weight $a=0$, the one descendant $\IZ_3$-Hodge
integrals in question are computed using the recursions of
\cite{cc:c3z3}. Integrating the Maple code\footnote{All codes can be made
  available to the interested reader upon request.} written by
Cadman-Cavalieri with formula (\ref{c3z3disc}) we recover all invariants in
Table 3.3 of \cite{Bouchard:2008gu} with the physical framing $f=-2/3$. \footnote{Further
computations (for $a= 2/3$), in agreement with the $f=0$ invariants of
\cite{Bouchard:2008gu} , suggest that the relationship should be a simple
translation $a=f+2/3$.}

\section{The main case: $[\bbC^3/\bbZ_4]$}
\label{sec:z4}
In this section we consider the orbifold $[\IC^3/\IZ_4]$, where the orbifold
group acts with weights $(1,1,2)$, for two different choices of
Lagrangians. When $\alpha_1=2$, $\alpha_0=\alpha_2=1$, the action is effective
along the axis that gets doubled to become the zero section of the
orbi-bundle. We refer to this choice of weights as the ``asymmetric
choice''. We then treat the case when $\alpha_0=2$,
$\alpha_1=\alpha_2=1$, in which the action is symmetric between the fibers and
has instead a non-trivial $\bbZ_2$ stabilizer along the base. 

\subsection{$B$-model, asymmetric case}
\label{sec:c3z4basymm}
\begin{figure}[h]
\begin{minipage}[t]{0.48\linewidth}
\centering
\vspace{0pt} 
\includegraphics[scale=0.5]{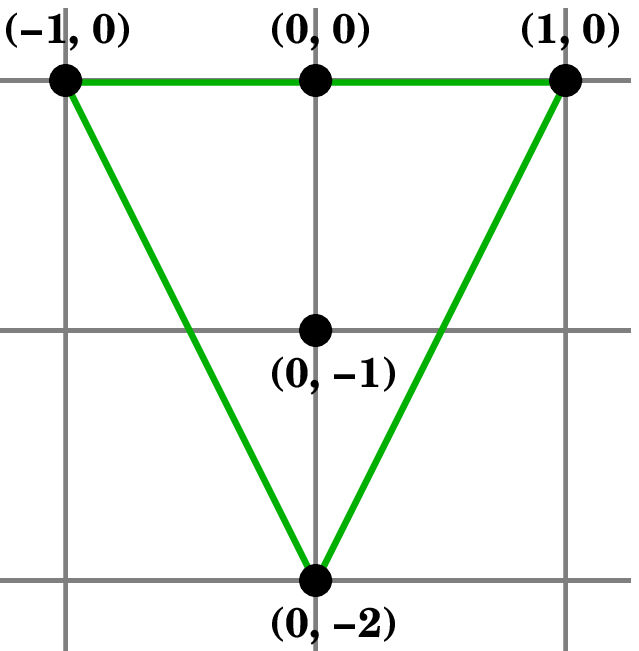}
\vspace{0pt} 
\caption{The Newton polytope of the fan of
  $\mathfrak{X}=[\bbC^3/\bbZ_4]$.}
\label{fig:c3z4}
\vspace{1.5cm}
\end{minipage} 
\begin{minipage}[t]{0.48\linewidth}
\centering
\vspace{0pt} 
\centering
  \includegraphics[scale=0.5]{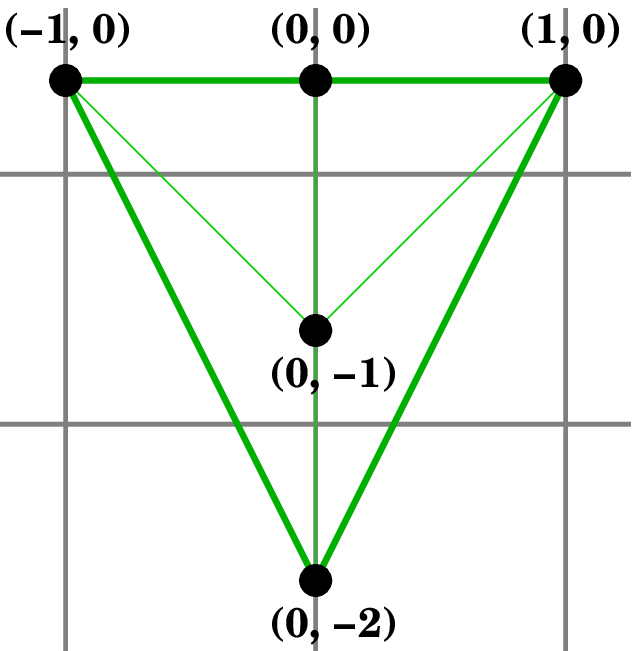}
\caption{The Newton polytope of the fan of $X=K_{\bbF_2}$.}
\label{fig:kf2}
\end{minipage}
\end{figure}

The Newton polytope associated to the fan of $[\bbC^3/\bbZ_4]$ is depicted in Fig. \ref{fig:c3z4}. Accordingly, the mirror curve has the following form
\beq
P_{[\bbC^3/\bbZ_4]}(U,V)=V+\frac{1}{V}-a_{\frac{1}{2}} - a_{\frac{1}{4}}/U - 1/U^2=0
\label{eq:mirrorcurvez4}
\eeq

In \eqref{eq:mirrorcurvez4}, $a_{\frac{1}{4}}$ and $a_{\frac{1}{2}}$ are the
$B$-model co-ordinates which are mirror to the small quantum co-homology
parameters $\tau_{\frac{1}{4}}$, $\tau_{\frac{1}{2}}$, where we write
$H^2_{\rm CR}([\bbC^3/\bbZ_4]) \ni \Phi = \tau_{\frac{1}{4}}
\mathbf{1}_{\frac{1}{4}}+\tau_{\frac{1}{2}} \mathbf{1}_{\frac{1}{2}}$. The
precise relation was found in \cite{Brini:2008rh, MR2486673}; we have
\bea
\tau_{\frac{1}{2}} &=& 2 \arcsin{\l( \frac{a_{\frac{1}{2}}}{2}\r) }
\eea
and at the first few orders in $a_{\frac{1}{4}}$, $a_{\frac{1}{2}}$
\bea
\tau_{\frac{1}{4}} &=&
\left(1+\frac{a_{\frac{1}{2}}^2}{32}+\frac{25 a_{\frac{1}{2}}^4}{6144}\right)
   a_{\frac{1}{4}}+\left(-\frac{a_{\frac{1}{2}}}{192}-\frac{25
   a_{\frac{1}{2}}^3}{18432}\right)
   a_{\frac{1}{4}}^3+ \dots \quad
\eea
\subsubsection{The $B$-model disc potential}
In writing \eqref{eq:mirrorcurvez4} we have implicitly made a choice of a
$SL(2,\mathbb{Z})$ representative for the spectral curve. It was argued in
\cite{Brini:2008rh} that this choice corresponds to the analytic continuation
at the orbifold point of an open string setup with branes on the upper legs of
the $pq$-web diagram of $K_{\bbF_2}$: this corresponds precisely to the asymmetric case for framing $f=1$.
 The Hori-Vafa differential $d\lambda_{[\bbC^3/\bbZ_4]}$ corresponding to \eqref{eq:mirrorcurvez4}, which gives the derivative of the $B$-model disc function \eqref{eq:bdisc}, reads
\bea
d\lambda^{(\alpha),
  f=1}_{[\bbC^3/\bbZ_4]}(a_{\frac{1}{4}},a_{\frac{1}{2}},U) &=& \log
\Bigg(\frac{1}{2U^4} \Bigg(1+a_{\frac{1}{4}} U+a_{\frac{1}{2}}U^2+ \nn \\ &+& \sqrt{\left(1+a_{\frac{1}{4}} U+a_{\frac{1}{2}}U^2\right)^2-4 U^4}\Bigg)\Bigg)\frac{dU}{U}
\label{eq:diffasymm}
\eea
whereas the open orbifold mirror map is trivial \cite{Brini:2008rh}
\beq
w=U
\eeq
up to the sign ambiguity that, as we have reviewed in Sec. \ref{sec:Bmodel},
is intrinsic in the definition of open invariants. We have appended a
superscript $(\alpha)$ to the differential to stress the fact that it refers to the
asymmetric choice.
The $B$-model disc potential is then 
\beq
\cF^{[\bbC^3/\bbZ_4], (\alpha, f=1)}_{0,1}(\tau_{\frac{1}{4}},\tau_{\frac{1}{2}},w)=\int^{U(w)} d\lambda^{(\alpha),f=1}_{[\bbC^3/\bbZ_4]}(a_{\frac{1}{4}}(\tau_{\frac{1}{4}}, \tau_{\frac{1}{2}}),a_{\frac{1}{2}}(\tau_{\frac{1}{2}}),U)
\label{eq:Bdiscpotasymm}
\eeq

\subsubsection{Higher genus open invariants from mirror symmetry}
\label{sec:higherB}
In this section we work out in detail the general predictions of open orbifold mirror symmetry for open invariants. This lays the basis for the comparision with the $A$-model computation of the orbifold annulus function in Sec.~\ref{sec:Aannulasymm}, and provides highly non-trivial predictions for some $g\leq 2$ open orbifold potentials. \\
We start with the following
\begin{thm}
Conjecture \ref{conj:modular} is true for $[\mathbb{C}^3/\bbZ_4]$. In this case $\mathcal{G}=\Gamma(2)$, i.e., the group of $SL(2,\bbZ)$ matrices congruent to the identity modulo 2.
\label{thm:c3z4mod}
\end{thm}
\noindent Some of the arguments to prove it were used, in a slightly different context, in \cite{Brini:2008rh}. We need the following technical
\begin{lemma}
\label{lem:WGprop}
  Let $\mathcal{S}=(\Gamma, C,\Delta, u,v)$ be a spectral curve with genus 1 support, i.e. $\mathfrak{g}(\Gamma)=1$, and logarithmic branch cuts $\Delta\neq \emptyset$. Let $V:=e^v$ be a degree 2 branched covering map to $\bbP^1$ and $q_1,q_2,q_3,q_4$ be its branch points; the fact that $V$ be degree 2 can always be accomplished up to a symplectic transformation \eqref{eq:gl2z}. Then the Eynard-Orantin correlators \eqref{eq:BKMPrecursion1} have the form for $3g-3+2h>0$
\beq
W_{h}^{(g)}(p_1 \ldots, p_h, \mathcal{S})=\sum_{l=0}^{3g-3+2h}A_l(p_1 \ldots, p_h, \{q_i\}_{i=1}^4) G^l(\{q_i\})
\label{eq:WGprop}
\eeq
where the {\rm propagator} $G(\{q_i\})$ is defined as
\beq
G(q_1,q_2,q_3,q_4)=\frac{E\l(k\r)}{K\l(k\r)}
\label{eq:Gprop}
\eeq
$$k=\frac{(q_1-q_3)(q_2-q_4)}{(q_1-q_2)(q_3-q_4)} $$

and $A_l(p_1 \ldots, p_h, \{q_i\})$ are for all $l$ meromorphic functions of $\{p_j\}$ for every $\{q_i\}$, and algebraic functions of the complex moduli of $v du$ for every $\{p_j\}$.
\end{lemma}
 The ordering of the set of branch points in \eqref{eq:Gprop} is dictated by the choice of polarization $\mathcal{H} \in H_1(\Gamma,\bbZ)$ of the spectral curve. In \eqref{eq:Gprop}, $E(x)$ and $K(x)$ denote the complete elliptic integrals of the second and first kind respectively.\\

{\noindent \it Proof.} We just sketch here the main lines of the proof; the interested reader may find the details in Appendix \ref{sec:appmodular}. \\
A proof can be given recursively. First of all \eqref{eq:WGprop} is true for the Bergmann kernel (as derived in \eqref{eq:bergmann}). Then, the Eynard-Orantin recursion straightforwardly implies \eqref{eq:WGprop} for $g=0$; when $g>0$, if we assume that \eqref{eq:WGprop} is true for $W_{h+1}^{(g-1)}(p, p_1 \ldots, p_h)$, then expressing the residues of the prime form $dE(p,q)$ in terms of elliptic integrals  shows that the expansion \eqref{eq:WGprop} holds for $W_{h}^{(g)}(p_1 \ldots, p_h)$ (see \eqref{eq:C}, \eqref{eq:Creg} and \eqref{eq:derpi}). By regularity of the curve, all coefficients $A(p_1 \ldots, p_h, \{q_i\})$ are algebraic in the complex moduli of $v(p) du(p)$; meromorphicity in $p_1 \ldots, p_h$ is trivially proven recursively.
\begin{flushright}$\square$\end{flushright}
    {\noindent \it Proof of Theorem \ref{thm:c3z4mod}.} From \eqref{eq:mirrorcurvez4}, we see that the family of curves for $[\mathbb{C}^3/\bbZ_4]$ is given as
\beq
\l(V+\frac{1}{V}\r)=a_{\frac{1}{2}}  + a_{\frac{1}{4}}/U + 1/U^2
\label{eq:torus}
\eeq
i.e., the support $\Gamma_X$ of the spectral curve is given by a family of
complex tori, and either $V$ or $U$ realize $\Gamma$ as a twofold branched
covering of $\bbP^1$. Therefore $\mathfrak{g}=1$; since $b_2^{\rm CR}([\bbC^3/\bbZ^4])=2$, we have one ``tilded'' period integral in the notation of Remark \ref{rem:noncpt} and Theorem \ref{thm:modular}, i.e., one flat co-ordinate of $QH^\bullet(K_{\bbF_2}) \simeq QH^\bullet([\bbC^3/\bbZ_4])$ which is not dual to a compact divisor. Closed mirror symmetry considerations \cite{Brini:2008rh} show that $\tilde t=a_{\frac{1}{2}}$. \\
By exploiting the analogy with the Seiberg-Witten curves of $\mathcal{N}=2$ pure Yang-Mills, it was shown in \cite{Brini:2008rh} that the branch points $q_i$ of the $V$-projection are given by
\beq
1/U=-\frac{a_{\frac{1}{4}}}{2}\pm c_1(\tau), \quad 1/U=-\frac{a_{\frac{1}{4}}}{2}\pm c_2(\tau)
\label{eq:bpts}
\eeq
where in terms of the elliptic modulus $\tau$ of the torus \eqref{eq:torus} we have
\bea
\label{eq:c1c2}
c_1(\tau)&=& \frac{2 \theta_4^2(\tau)}{ \theta_2^2(\tau)} \\
 c_2(\tau) &=& \frac{2 \theta_3^2(\tau)}{ \theta_2^2(\tau)} \\
 a_{\frac{1}{4}}(a_{\frac{1}{2}},\tau) &=& 2\sqrt{\frac{4\theta_4^2(\tau)}{ \theta_2^2(\tau)}+2+a_{\frac{1}{2}}}
\label{eq:aimod}
\eea
By Lemma \ref{lem:WGprop} and \eqref{eq:c1c2}-\eqref{eq:aimod} we have
\beq
W_{h}^{(g)}(p_1 \ldots, p_h)=\sum_{l=0}^{3g-3+2h}A_l(p_1 \ldots, p_h, a_{\frac{1}{2}}, \tau) G^l(c_1,c_2)
\eeq
where we have expressed the closed modulus $a_{\frac{1}{4}}$ as a function of the non-compact period $a_{\frac{1}{2}}$ and the elliptic modulus, and the dependence on $a_{\frac{1}{2}}$ cancels from the propagator because of \eqref{eq:bpts}. \\
$A(p_1 \ldots, p_h, a_{\frac{1}{2}}, \tau)$ is then, for every $\{p_i\}$ and $a_{\frac{1}{2}}$, a holomorphic weight zero modular form of $\Gamma(2)$ by \eqref{eq:c1c2}-\eqref{eq:aimod}, since the Jacobi theta functions are modular forms of $\Gamma(2)$ of weight $1/2$. As far as $G(c_1,c_2)$ is concerned we use the fact that, denoting 
$$k=\frac{(q_1-q_3)(q_2-q_4)}{(q_1-q_2)(q_3-q_4)} $$
we have the remarkable identities
\beq
E(k) K(k)=\l(\frac{\pi}{2}\r)\frac{4E_2(2\tau)-E_2(\tau)}{3}
\eeq
and
\beq
K(k)= \frac{2}{\pi}\theta_4^2(\tau)
\eeq
Using the duplication formula
\beq
E_2(2\tau)=\frac{2E_2(\tau)+\theta_3^2(\tau)+\theta_4^2(\tau)}{4}
\eeq
the claim follows.
\begin{flushright}$\square$\end{flushright}
We  now turn to the $B$-model computation of higher order open orbifold
potentials. To this aim, let us fix first a choice of polarization
$\mathcal{H}\in H_1(\Gamma_{[\bbC^3/\bbZ^4]},\bbZ)$: the point of maximally
unipotent monodromy \cite{MR1677117} of the torically compactified $B$-model
moduli space is given, in inhomogeneous $B$-model co-ordinates, by
$(a_{\frac{1}{2}}, a_{\frac{1}{4}})\sim(\infty,\infty)$
\cite{Brini:2008rh}. We will fix a polarization of $\mathcal{H}$ as follows:
let $A$ (resp. $B$) $\in H_1(\Gamma_{[\bbC^3/\bbZ^4]}, \bbZ)$ be the 1-cycle represented by a
loop encircling the $[q_1, q_2]$ (resp. $[q_3, q_4]$) segment in the
$U$-plane. We order the set of branch points $\{q_i\}_{i=1}^4$ such that the
periods of $d\lambda_{[\bbC^3/\bbZ_4}]$ around $A$ (resp. $B$) has a logarithmic
(resp. double-logarithmic) singularity around the maximally unipotent monodromy point. This corresponds to computing $W_{h}^{(g)}$ in the so-called ``large radius phase''. Then we have the following
\begin{prop}
Let $W_{h}^{(g)}(p_1 \ldots, p_h)$ be the correlators computed from the recursion \eqref{eq:BKMPrecursion1}-\eqref{eq:BKMPrecursion3} with the choice of polarization $\mathcal{H}$ above, and let \eqref{eq:WGprop} be their polynomial expansion in powers of the propagator. Then the {\rm orbifold correlators} $\mathcal{W}_{h}^{(g)}(p_1 \ldots, p_h)$ are given by
\beq
\mathcal{W}_{h}^{(g)}(p_1 \ldots, p_h)=\sum_{l=0}^{3g-3+2h}A_l(p_1 \ldots, p_h, \{q_i\}) \tilde{G}^l(\{q_i\})
\label{eq:WGproporb}
\eeq
where the coefficients $A_l(p_1 \ldots, p_h, \{q_i\})$ coincide with those in \eqref{eq:WGprop}, and the {\rm orbifold propagator} $\tilde G$ is defined by
\beq
\label{eq:orbprop}
\tilde G(a_{\frac{1}{4}},a_{\frac{1}{2}}) := G(a_{\frac{1}{4}},a_{\frac{1}{2}})
-\frac{\pi }{2 K\left(w^-\right)
   \left(K\left(w^-\right)+K\left(w^+\right)\right)}
\eeq
where
\bea
w^- &=& -\frac{\left(\sqrt{a_{\frac{1}{4}}^2-4
   a_{\frac{1}{2}}-8}-\sqrt{a_{\frac{1}{4}}^2-4
   a_{\frac{1}{2}}+8}\right)^2}{4 \sqrt{\left(a_{\frac{1}{4}}^2-4
   a_{\frac{1}{2}}\right)^2-64}}, \\
w^+ &=& \frac{\left(\sqrt{a_{\frac{1}{4}}^2-4 a_{\frac{1}{2}}-8}+\sqrt{a_{\frac{1}{4}}^2-4
   a_{\frac{1}{2}}+8}\right)^2}{4 \sqrt{\left(a_{\frac{1}{4}}^2-4
   a_{\frac{1}{2}}\right)^2-64}}
\eea
\end{prop}
The proof relies on applying the transformation \eqref{eq:WhgE2trasf} with the shift \eqref{eq:dtau} and the change of basis \eqref{eq:M2} computed in \cite{Brini:2008rh}

\beq
M=
\left(
\begin{array}{cc}
\frac{2 \pi ^{3/2}}{\Gamma \left(\frac{1}{4}\right)^2}
&
\frac{(1-i) \sqrt{\pi }}{\Gamma \left(\frac{1}{4}\right)^2} \\
-\frac{\Gamma \left(\frac{1}{4}\right)^2}{\sqrt{\pi }}
&
\frac{\left(\frac{1}{2}+\frac{i}{2}\right) \Gamma \left(\frac{1}{4}\right)^2}{\pi ^{3/2}}
\end{array}
\right)
\label{eq:MZ4}
\eeq
It is quite remarkable to notice that the sole analytic continuation of the ``large radius'' open string generating functions around $(a_{\frac{1}{2}},a_{\frac{1}{4}})\simeq (0,0)$, without the shift $E_2(\tau) \to E_2(\tau)+d(\tau)$, would end up in an expansion in $\tau_{\frac{1}{2}},\tau_{\frac{1}{4}}$ with irrational (in fact transcendental) coefficients. Indeed,  the propagator $G$ has the following expansion in flat co-ordinates
\bea
G(\tau_{\frac{1}{4}},\tau_{\frac{1}{2}}) &=& 
\left(\frac{1}{2}+\frac{4 \pi ^2}{\Gamma \left(\frac{1}{4}\right)^4}\right)+\left(-\frac{i}{32}-\frac{2 i \pi ^4}{\Gamma
   \left(\frac{1}{4}\right)^8}\right) \tau_{\frac{1}{4}}^2  \\ &+& \left(\left(\frac{i}{8}+\frac{8 i \pi
   ^4}{\Gamma \left(\frac{1}{4}\right)^8}\right)+\left(\frac{8 \pi ^6}{\Gamma \left(\frac{1}{4}\right)^{12}}-\frac{\pi ^2}{8
   \Gamma \left(\frac{1}{4}\right)^4}\right) \tau_{\frac{1}{4}}^2\right) \tau_{\frac{1}{2}} \nn \\ &+& \left(\left(-\frac{16 \pi
   ^6}{\Gamma \left(\frac{1}{4}\right)^{12}}+\frac{\pi ^2}{4 \Gamma
   \left(\frac{1}{4}\right)^4}\right)+
\right)
   \tau_{\frac{1}{2}}^2+ \dots \nn
\eea
The terms containing powers of $\Gamma(1/4)$ are exactly cancelled by the shift in the propagator in \eqref{eq:orbprop}
\bea
\tilde G(\tau_{\frac{1}{4}}, \tau_{\frac{1}{2}})  &=& \frac{1}{2}-\frac{i \tau_{\frac{1}{4}}^2}{32}-\frac{11 i
   \tau_{\frac{1}{4}}^6}{61440}+\left(
   \frac{i}{8}+\frac{13 i
   \tau_{\frac{1}{4}}^4}{6144}+\frac{457 i
   \tau_{\frac{1}{4}}^8}{13762560}\right)
   \tau_{\frac{1}{2}} \nn \\ &+& \left(-\frac{i
   \tau_{\frac{1}{4}}^2}{128}-\frac{371 i
   \tau_{\frac{1}{4}}^6}{1474560}\right)
   \tau_{\frac{1}{2}}^2 + \dots
\eea
Conjecture \ref{conj:BKMP2} then implies that, upon integrating with respect
to $p_1,\dots, p_h$ and plugging in the mirror map, the orbifold correlators
should provide the genus $g$, $h$-holes orbifold potentials of
$[\bbC^3/\bbZ_4]$. The results of our $B$-model computations of open orbifold
Gromov-Witten invariants up to genus 2 are contained in Appendix
\ref{sec:numbers}.

\subsection{$A$-model, asymmetric case}
\label{sec:c3z4as}
 It appears that the way to compare the localization computations with
 the $B$-model predictions at $f=1$ is to choose $a$ to be equal to one over the effective degree of the action in the first fiber direction. In this case the torus weights become:
\beq
\label{eq:weightsasymm}
(s_0,s_1,s_2)=\left(\frac{1}{4}, -\frac{1}{2}, \frac{1}{4}\right).
\eeq
 
Insertions that give rise to non-equivariant invariants correspond to the two
age one orbifold cohomology classes, $\mathbf{1}_{\frac{1}{4}},
\mathbf{1}_{\frac{2}{4}}$. The compatibility condition between degree and twisting (\ref{degtwist}) is $k\equiv d$ (mod $4$), and the disc function is:
\beq
\label{eq:Adiscfunasymm}
\mathcal{D}^{(\alpha)}(d,1/2)= \frac{1}{\left\lfloor \frac{d}{4}\right\rfloor!}\left(\frac{\hbar}{d}\right)^{\mathrm{age}(\mathbf{1}_{\frac{d}{4}})}\frac{\Gamma(\frac{d}{4}+\langle \frac{d}{4}\rangle-\frac{d}{2})}{\Gamma(1-\langle \frac{d}{4}\rangle-\frac{d}{2})} (\mathbf{1}_{\frac{d}{4}})^\vee
\eeq
 Once again, we have added a superscript $(\alpha)$ to stress the fact that we
 refer to the asymmetric choice.

\subsubsection{A mirror theorem for orbifold disc invariants}
The small $J$- function for the closed theory is:

\beq
\begin{aligned}
J(s_i; \tau_{\alpha};z) = z\mathbf{1}+\tau_{\frac{1}{4}}\mathbf{1}_{\frac{1}{4}}+\tau_{\frac{2}{4}}\mathbf{1}_{\frac{2}{4}}+
  \sum_{m_1,
    m_2}\sum_{k=0}^3\frac{\tau_{\frac{1}{4}}^{m_1}\tau_{\frac{2}{4}}^{m_2}(\mathbf{1}_{\frac{k}{4}})^\vee}{m_1!m_2!}
  \\
\times \int_{\smbzf{m_1+m_2+1}{\mathbf{1}_{\frac{1}{4}}^{m_1},\mathbf{1}_{\frac{2}{4}}^{m_2}, \mathbf{1}_{\frac{k}{4}}}}\hspace{-1.5cm} \frac{e(\bE_{1}^\vee\otimes\mathcal{O}(s_0)\oplus\bE_{2}^\vee\otimes\mathcal{O}(s_1)\oplus \bE_{1}^\vee\otimes\mathcal{O}(s_2))}{(z - \psi)}
\end{aligned}
\eeq

and the potential for open disc invariants is given by:

\beq
\begin{aligned}
F_{0,1}^{[\bbC^3/\bbZ_4], (\alpha)}(\tau_\frac{1}{4},\tau_\frac{2}{4},w,{1}/{2}) := \sum_{m_1,m_2,d}\langle\mathbf{1}_{\frac{1}{4}}^{m_1}\mathbf{1}_{\frac{2}{4}}^{m_2};\frac{1}{2}\rangle^{d}_{0}\frac{\tau_{\frac{1}{4}}^{m_1}}{{m_1}!}\frac{\tau_{\frac{2}{4}}^{m_2}}{m_2!}\frac{w^d}{d!} \nonumber\\
= \sum_d\left[J\left(\frac{1}{4},-\frac{1}{2},\frac{1}{4};0,\tau_\frac{1}{4},\tau_\frac{2}{4};\frac{1}{d}\right)\mathcal{D}^{(\alpha)}(d,{1}/{2}) \right]\frac{w^d}{d!}
\label{eq:Adiscpotasymm}
\end{aligned}
\eeq

As was argued in \cite{Brini:2008rh} by physical considerations of monodromy invariance, the asymmetric case is the one for which the $B$-model computations of the previous section have the best chance to yield a correct answer. In this section we give a full proof of this statement, by establishing a version of open orbifold mirror symmetry for disc invariants for this example.
\begin{thm}
The analytic part of the $B$-model disc potential \eqref{eq:Bdiscpotasymm} coincides, for positive winding numbers and up to signs, with the generating function of orbifold Gromov-Witten disc invariants \eqref{eq:Adiscpotasymm}.
\label{thm:mirror}
\end{thm}
\begin{rem}
Theorem \ref{thm:mirror} postulates open orbifold mirror symmetry as an almost-equality of $A$ and $B$-model open string potentials, which coincide only after dropping non-analytic terms and up to signs. This might seem a bit of a nuisance, but it should in fact be entirely expected: the non-analytic terms that are dropped, analogous to the power-of-a-logarithm terms of their closed string counterparts, are degree-zero contributions for which we do not have a clear $A$-model definition, and likewise for the zero-winding number term. On the other hand, the possible sign differences reside in the inherent ambiguity in the definition of the open string potential and mirror map on the $B$-model side.
In this case, again, the unfixed torus weight $a$ is identified with the framing ambiguity on the mirror side.
\end{rem}
{\noindent \it Proof of Theorem \ref{thm:mirror}.} We explicitly evaluate the power series expansion of the $A$-model disc function \eqref{eq:Adiscpotasymm} in the winding parameter by analyzing the expression of the twisted equivariant $J$-function of $[\bbC^3/\bbZ_4]$, and compare the results with the analogous expansion of the $B$-model disc function as written in \eqref{eq:Bdiscpotasymm}. The key idea is to work with closed $B$-model co-ordinates, i.e. $a_{\frac{1}{4}}$ and $a_{\frac{1}{2}}$, instead of flat ones. To begin with, define 
\bea g^A(a_{\frac{1}{4}},a_{\frac{1}{2}},w) &:=& \de_w \de_{a_{\frac{1}{4}}}F^{[\bbC^3/\bbZ_4], (\alpha)}_{0,1}(\tau_{\frac{1}{4}}(a_{\frac{1}{4}}, a_{\frac{1}{2}}),\tau_{\frac{1}{2}}(a_{\frac{1}{2}}),w) \\ g^B(a_{\frac{1}{4}},a_{\frac{1}{2}},w) &:=& \de_w \de_{a_{\frac{1}{4}}}\cF^{[\bbC^3/\bbZ_4], (\alpha)}_{0,1}(\tau_{\frac{1}{4}}(a_{\frac{1}{4}}, a_{\frac{1}{2}}),\tau_{\frac{1}{2}}(a_{\frac{1}{2}}),w)
\eea
From \eqref{eq:diffasymm} and \eqref{eq:Bdiscpotasymm}, we obtain for $g^B$
\beq
g^B(a_{\frac{1}{4}}, a_{\frac{1}{2}},w)=\frac{1}{\sqrt{\left(1+a_{\frac{1}{4}} w+a_{\frac{1}{2}}w^2\right)^2-4 w^4}}
\eeq
The expression above simplifies greatly the task of finding a closed form for the Taylor coefficients of $\cF^{[\bbC^3/\bbZ_4], (\alpha)}_{0,1}$. Expansion of the square root at the denominator around $w=0$, Newton's binomial formula and standard power series manipulations yield
\bea
g^B(a_{\frac{1}{4}}, a_{\frac{1}{2}},w) &=& \sum_{m,n,k=0}^\infty (-1)^k 4^n \left(
\begin{array}{c}
 k \\
 -m+2 k+4 n
\end{array}
\right) \left(\frac{1}{2}\right)_n (2 n+1)_k \nn \\ & \times &  \frac{a_{\frac{1}{2}}^{-m+2 k+4 n} a_{\frac{1}{4}}^{m-k-4 n}}{k! n!}
\label{eq:gB}
\eea
where as usual the binomial function and the Pochhammer symbol are defined as
$$\l( \begin{array}{c}
 n \\
k
\end{array}\r) = \frac{n!}{k!(n-k)!}, \quad (a)_n= \frac{\Gamma(n+a)}{\Gamma(a)} $$
Let us turn to analyze Eq. \eqref{eq:Adiscpotasymm}. 
The small $J$-function of $[\bbC^3/\bbZ_4]$ in $B$-model co-ordinates is given \cite{MR2276766, MR2510741, MR2486673} by the following expression
\beq
J\l(a_{\frac{1}{4}},a_{\frac{1}{2}};z\r)=z \sum_{m,n=0}^\infty a_{\frac{1}{4}}^m a_{\frac{1}{2}}^n R_{-n/2-m/4,-m/2}(z) \mathbf{1}_{\langle n/2+m/4 \rangle}
\label{eq:Jfunasymm}
\eeq
where for the asymmetric case and weights \eqref{eq:weightsasymm} we have
\beq
R_{k,l}(z)= \frac{\prod_{b=k+1}^0(\frac{1}{4}+ b z)^2 \prod_{b=l+1}^0(-\frac{1}{2} +b z)}{(l-2k)!(-2 l)! z^{-2k-l}}
\label{eq:RJfun}
\eeq
and we have denoted with $\langle x\rangle=x-[x]$ the fractional part of a real
number $x$ and with the short-hand notation $J\l(a_{\frac{1}{4}},a_{\frac{1}{2}};z\r):=J\l(\frac{1}{4},-\frac{1}{2},\frac{1}{4};0;\tau_{\frac{1}{4}}(a_{\frac{1}{4}},a_{\frac{1}{2}}), \tau_{\frac{1}{2}}(a_{\frac{1}{2}});z\r)$ the small $J$-function, expressed in $B$-model co-ordinates, with the torus weights given by \eqref{eq:weightsasymm}. \\
By the form \eqref{eq:Adiscfunasymm} of the disc function, the contribution of the $J$-function to the winding number $d$ term of $F^{[\bbC^3/\bbZ_4], (\alpha,f)}_{0,1}$ comes from the component of $J$ proportional to $\mathbf{1}_{\frac{d}{4}}$. It is therefore convenient to isolate the projection of the $J$-function to $\mathbf{1}_{\frac{d}{4}}$ at winding number $d=k \hbox{ mod } 4$ for each $k=0,1,2,3$. To this aim, denote the projections as
$$J=:\sum_{k=0}^3 J_{[k]} \mathbf{1}_{\frac{k}{d}} $$
For $k=0$, i.e. $d=4 L$, $L\in \bbZ$, we find from \eqref{eq:Jfunasymm}, \eqref{eq:RJfun}
\bea
J_{[0]} &=: &J_{[0]}^{\mathrm{even}}+  J_{[0]}^{\mathrm{odd}} \nn \\
 J_{[0]}^{\mathrm{even}}\l(a_{\frac{1}{4}},a_{\frac{1}{2}},\frac{1}{4L}\r) &=& \sum_{a,b=0}^\infty\Bigg[\frac{\Gamma (-2 L) \Gamma (L+1)^2 a_{\frac{1}{4}}^{4 b+2} a_{\frac{1}{2}}^{2 a+1}}{2 (2
   a+1)! (4 b+2)! \Gamma (-a-b+L)^2} \nn \\ & & \frac{1}{\Gamma (-2 (b+L))}\Bigg] \\
 J_{[0]}^{\mathrm{odd}}\l(a_{\frac{1}{4}},a_{\frac{1}{2}},\frac{1}{4L}\r) &=& \sum_{a,b=0}^\infty \Bigg[\frac{\Gamma (1-2 L) \Gamma (L+1)^2 a_{\frac{1}{4}}^{4 b} a_{\frac{1}{2}}^{2 a}}{4 L (2
   a)! (4 b)! \Gamma (-2 b-2 L+1)} \nn \\ & & \frac{1}{ \Gamma (-a-b+L+1)^2}\Bigg]
\eea
while the disc function is
\beq
\mathcal{D}^{(\alpha)}(4L,1/2)= \frac{(-1)^L (2 L-1)!}{(L!)^2} \l(\mathbf{1}_0\r)^\vee
\eeq
For $k=4L+1$ we obtain likewise
\bea
J_{[1]} &=:& J_{[1]}^{\mathrm{even}}+  J_{[1]}^{\mathrm{odd}} \nn \\
\label{eq:Jfirst}
J_{[1]}^{\mathrm{even}}\l(a_{\frac{1}{4}},a_{\frac{1}{2}},\frac{1}{4L+1}\r) &=& \sum_{a,b=0}^\infty\Bigg[
\frac{(-2 (b+L))_{2 b}}{(2 a)! \Gamma (4 b+2)}
a_{\frac{1}{4}}^{4 b+1} a_{\frac{1}{2}}^{2 a} \nn \\ & & (-a-b+L+1)_{a+b}{}^2 \Bigg] \\
J_{[1]}^{\mathrm{odd}}\l(a_{\frac{1}{4}},a_{\frac{1}{2}},\frac{1}{4L+1}\r) &=& \sum_{a,b=0}^\infty\Bigg[
\frac{ (-a-b+L)_{a+b+1}{}^2}{\Gamma (2 a+2) \Gamma (4
   b+4)}
 a_{\frac{1}{4}}^{4 b+3} a_{\frac{1}{2}}^{2 a+1} \nn \\ & & (-2 b-2 L-1)_{2 b+1} \Bigg] \nn \\
\eea
and
\beq
\mathcal{D}^{(\alpha)}(4L+1,1/2)=\frac{2 (-1)^L \Gamma (2 L)}{(4 L+1) L! \Gamma (L)}\l(\mathbf{1}_{\frac{1}{4}}\r)^\vee
\eeq
For $k=4L+2$ we have in the same way
\bea
J_{[2]} &=:& J_{[2]}^{\mathrm{even}}+  J_{[2]}^{\mathrm{odd}} \nn \\
J_{[2]}^{\mathrm{even}}\l(a_{\frac{1}{4}},a_{\frac{1}{2}},\frac{1}{4L+2}\r) &=& \sum_{a,b=0}^\infty\Bigg[
\frac{  \left(-2 b-\frac{1}{8
   L+4}\right)_{2 b+1}}{(2 a)! \Gamma (4 b+3)} a_{\frac{1}{2}}^{2a} a_{\frac{1}{4}}^{4b+2}
\nn \\ & & \left(- a- b+\frac{L+1}{4
   L+2}\right)_{a+b}^2 \Bigg] \\
J_{[2]}^{\mathrm{odd}}\l(a_{\frac{1}{4}},a_{\frac{1}{2}},\frac{1}{4L+2}\r) &=& \sum_{a,b=0}^\infty\Bigg[
\frac{ \left(-2 b-\frac{1}{8
   L+4}+1\right)_{2 b}}{(4 b)! \Gamma (2 a+2)} a_{\frac{1}{2}}^{2a+1}
a_{\frac{1}{4}}^{4b} \nn \\ & & \left(- a- b+\frac{L+1}{4
      L+2}\right)_{a+b}^2\Bigg]
\eea
and
\beq
\mathcal{D}^{(\alpha)}(4L+2,1/2)=\frac{2 (-1)^L \Gamma (2 L)}{(4 L+2) L! \Gamma (L)}\l(\mathbf{1}_{\frac{1}{2}}\r)^\vee
\eeq
Finally for $k=4L+3$
\bea
J_{[3]} &=:& J_{[3]}^{\mathrm{even}}+  J_{[3]}^{\mathrm{odd}} \nn \\
J_{[3]}^{\mathrm{even}}\l(a_{\frac{1}{4}},a_{\frac{1}{2}},\frac{1}{4L+2}\r) &=& \sum_{a,b=0}^\infty\Bigg[
\frac{  \left(\frac{2 L+1}{4
   L+3}-2 b\right)_{2 b}a_{\frac{1}{2}}^{2a+1} a_{\frac{1}{4}}^{4b+1}}{(4 L+3) \Gamma (2 a+2) \Gamma
   (4 b+2)}  
\nn \\ & & \left(- (a+ b)+\frac{L+1}{4
   L+3}\right)_{a+b}^2\Bigg]\\
J_{[3]}^{\mathrm{odd}}\l(a_{\frac{1}{4}},a_{\frac{1}{2}},\frac{1}{4L+2}\r) &=& \sum_{a,b=0}^\infty \Bigg[
\frac{  \left(\frac{2 (L+b (4
   L+3)+1)}{-4 L-3}\right)_{2 b+1}}{(4 L+3) (2 a)!
   \Gamma (4 b+4)}  a_{\frac{1}{2}}^{2a} a_{\frac{1}{4}}^{4b+3} \nn \\ & & \left(- a- b+\frac{L+1}{4
   L+3}\right)_{a+b}^2\Bigg]
\label{eq:Jlast}
\eea
and
\beq
\mathcal{D}^{(\alpha)}(4L+3,1/2)=\frac{(-1)^{L+1} \Gamma (2 L+2)}{L (4 L+3)^2 L! \Gamma
   (L)}\l(\mathbf{1}_{\frac{3}{4}}\r)^\vee
\eeq
With these expression at hand we can now make a detailed comparison with the $B$-model disc function. Write
\beq
g^A(a_{\frac{1}{4}},a_{\frac{1}{2}},w)=:\sum_{d=0}^\infty g^A_d(a_{\frac{1}{4}},a_{\frac{1}{2}}) w^d, \quad g^B(a_{\frac{1}{4}},a_{\frac{1}{2}},w)=:\sum_{d=0}^\infty g^B_d(a_{\frac{1}{4}},a_{\frac{1}{2}}) w^d
\eeq
We find from \eqref{eq:gB} and \eqref{eq:Jfirst}-\eqref{eq:Jlast}
\beq
g^A_{4L+k}(a_{\frac{1}{4}},a_{\frac{1}{2}})=(-1)^L g^B_{4L+k}(a_{\frac{1}{4}},a_{\frac{1}{2}}), \quad k=0,1,2,3
\eeq
This establishes mirror symmetry for all disc invariants with at least one insertion of $\mathbf{1}_{1/4}$. \\ As far as invariants with only $\mathbf{1}_{1/2}$-insertions are concerned, define $$h^A(a_{\frac{1}{2}},w):=\de_w \de_{a_{\frac{1}{2}}}F^{[\bbC^3/\bbZ_4], (\alpha)}_{0,1}(0,a_{\frac{1}{2}},w)$$ $$h^B(a_{\frac{1}{2}},w):=\de_w \de_{a_{\frac{1}{2}}}\cF^{[\bbC^3/\bbZ_4], (\alpha)}_{0,1}(0,a_{\frac{1}{2}},w)$$
On the $B$-model side the situation parallels closely what we have already done , given that
\beq 
\de_w \de_{a_{\frac{1}{2}}}\cF^{[\bbC^3/\bbZ_4], (\alpha)}_{0,1}(a_{\frac{1}{4}},a_{\frac{1}{2}},w)=w \de_w \de_{a_{\frac{1}{4}}}\cF^{[\bbC^3/\bbZ_4], (\alpha)}_{0,1}(a_{\frac{1}{4}},a_{\frac{1}{2}},w)
\eeq
as the reader can easily check, while on the $A$-model side we just have to isolate the $\mathcal{O}(a_{\frac{1}{4}}^0)$ terms in the $J$-function in order to compute $h^A$. We find, defining $h^A(a_{\frac{1}{2}},w)=:\sum_{d=0}^\infty h^A_d(a_{\frac{1}{2}}) w^d$, $h^B(a_{\frac{1}{2}},w)=:\sum_{d=0}^\infty h^B_d(a_{\frac{1}{2}}) w^d$, that
\beq
h^A_{4L+k}(a_{\frac{1}{4}},a_{\frac{1}{2}})=(-1)^L h^B_{4L+k}(a_{\frac{1}{4}},a_{\frac{1}{2}}), \quad k=0,2
\eeq
where the identity above is trivially true for $k=1,3$, when both sides are zero. \\
All we are left to do to complete the proof is to compute the $A$-model disc function in absence of insertion of twisted classes. Putting $a_{\frac{1}{4}}=0$, $a_{\frac{1}{2}}=0$ in \eqref{eq:Jfirst}-\eqref{eq:Jlast} and performing the sum over winding numbers in \eqref{eq:Adiscpotasymm} we find
\bea
\de_w F^{[\bbC^3/\bbZ_4], (\alpha)}_{0,1}(0,0,w) &=& \sum_{d=1}^\infty
(-1)^{d+1} \frac{(2d-1)!}{(d!)^2}w^{4d-1} \nn \\
&=&\frac{\log \left(\frac{1}{2} \left(\sqrt{4 w^4+1}+1\right)\right)}{w}
\eea
which, by \eqref{eq:diffasymm} and \eqref{eq:Bdiscpotasymm}, coincides with $\de_w F^{[\bbC^3/\bbZ_4], (\alpha)}_{0,1}(0,0,-w)$, upon dropping the non-analytic logarithmic term (i.e., restricting to positive degrees only). This concludes the proof.
\begin{flushright}$\square$\end{flushright}

\subsubsection{Annulus Invariants}
\label{sec:Aannulasymm}
Our localization formula expresses annulus invariants in terms of the disc function and of compact invariants with two descendant insertions:

\begin{equation}
\begin{array}{l}
{\langle \mathbf{1}^{m_0} \mathbf{1}_{\frac{1}{4}}^{m_1} \mathbf{1}_{\frac{2}{4}}^{m_{2}}\rangle_0^{d_1,d_2}=}  \\ \\
{\left(\frac{4}{\hbar}\right)^2 \prod_{j=1}^2 \left(\varphi(d_j) \DD^{(\alpha)}_{d_j}\left(d_j, 1/2\right)\right)
\int_{\mathcal{M}}\frac{e^{\mathrm{eq}}(\IE^\vee_{1}(1/4)\oplus\IE^\vee_{2}(-1/2)\oplus\IE^\vee_{1}(1/4))}{ \left({\frac{\hbar}{d_1}}-\psi_1\right)\left({\frac{\hbar}{d_2}}-\psi_2\right)}},
\end{array} 
\end{equation}
where 
$$
\mathcal{M}=\overline{\mathcal{M}}_{0,\sum m_j+2}(B\IZ_4,0;\mathbf{1}^{m_0} \mathbf{1}_{\frac{1}{4}}^{m_1} \mathbf{1}_{\frac{2}{4}}^{m_{2}},\mathbf{1}_{\frac{4-k_1}{4}}, \mathbf{1}_{\frac{4-k_2}{4}})
$$
and
$$
\varphi(k):=
\left\{
\begin{array}{cl}
-\frac{\hbar^3}{32} & k=0 \\
1                  & k=1,3 \\
\frac{\hbar}{4}    &k=2
\end{array}
\right.
$$

The second descendant insertion can be removed inductively using the genus $0$
topological recursion relations; this allows us to compute many invariants,
for which we find perfect
agreement with the mirror symmetry predictions of Sec. \ref{sec:higherB}. In Table \ref{tab:annulus1}-\ref{tab:annulus3} we collect the first
values for up to $n<9$ insertions. \\
We conclude this section with a very explicit example, to point out how to unravel explicitly the localization formula.

\begin{example}
We compute the annulus invariant
$\langle \mathbf{1}_\frac{1}{4}^2 \rangle_0^{(1,1)}$:
$$ \langle \mathbf{1}_\frac{1}{4} ^2 \rangle_0^{(1,1)}= 4^2 \int_{\mathbf{1}_\frac{1}{4}^2\mathbf{1}_\frac{3}{4}^2}\frac{\mathbb{E}_1^\vee(1/4)\oplus\mathbb{E}_{2}^\vee(-1/2)\oplus\mathbb{E}_1^\vee(1/4)}{(\hbar-\psi_3)(\hbar-\psi_4)} =$$
$$
\left[-\frac{1}{2}(\psi_3+\psi_4)+ 4(\lambda_1)_1-(\lambda_1)_{2}  \right]_{\mathbf{1}_\frac{1}{4}^2\mathbf{1}_\frac{3}{4}^2}=-\frac{1}{8},
$$
where the final evaluation is obtained via the explicit Hodge integral computations:
\begin{enumerate}
	\item $(\lambda_1)_{1}= \frac{1}{16}$;
	\item $(\lambda_1)_{2}=\frac{1}{8}$;
	\item $\psi_3=\psi_4=\frac{1}{4}$.
\end{enumerate}

\end{example}

\begin{table}[h]
\centering
\begin{tabular}{|c|ccccc|}
\hline
 & $m$ & 0 & 2 & 4 & 6  \\
\hline
$n$ & &  & & & \\
0 & & 0 & $-\frac{1}{8}$ & 0 & $-\frac{3}{128}$ \\
1 & & $\frac{1}{4}$ & 0 & $\frac{1}{128}$ & 0 \\
2 & & 0 & 0 & 0 & -$\frac{35}{512}$ \\ 
3 & & -$\frac{1}{32}$ & 0 & $\frac{3}{128}$ & 0 \\ 
4 & & 0 & $-\frac{3}{256}$ & 0 &   \\ 
5 & & $\frac{1}{64}$ & 0 &  &  \\ 
6 & & 0 &  &  & 
 \\
7 & & -$\frac{23}{1024}$ &  
&  &  \\
& &  & & & \\
\hline
\end{tabular}
\caption{Annulus orbifold Gromov--Witten invariants $\langle
  \mathbf{1}_\frac{1}{4}^m   \mathbf{1}_\frac{1}{2}^n \rangle_0^{(1,1)}$ of
  $[\bbC^3/\bbZ_4]$  in the asymmetric case at winding number $(1,1)$.}
\label{tab:annulus1}
\end{table}

\begin{table}[h]
\centering
\begin{tabular}{|c|ccccc|}
\hline
 & $m$ & 1 & 3 & 5 & 7 \\
\hline
$n$ & &  & & & \\
0 & & 0 & $\frac{3}{16}$ & 0 & $\frac{21}{1024}$  \\
1 & & -$\frac{1}{8}$ & 0 & $\frac{3}{512}$ & 0  \\
2 & & 0 & -$\frac{5}{256}$ & 0 &   \\
3 & & $\frac{5}{128}$ & 0 & -$\frac{255}{8192}$ &  \\
4 & & 0 & $\frac{35}{4096}$ & 0 &   \\
5 & & -$\frac{21}{2048}$ & 0 &  &   \\
6 & & 0 &  &  &   \\
7 & & $\frac{465}{32768}$ &  &  &   \\
& &  & & & \\
\hline
\end{tabular}
\caption{Annulus orbifold Gromov--Witten invariants $\langle
  \mathbf{1}_\frac{1}{4}^m   \mathbf{1}_\frac{1}{2}^n \rangle_0^{(2,1)}$of
  $[\bbC^3/\bbZ_4]$  in the asymmetric case at winding number $(2,1)$.}
\end{table}

\begin{table}[h]
\centering
\begin{tabular}{|c|ccccc|}
\hline
 & $m$ & 0 & 2 & 4 & 6 \\  
\hline
$n$ & &  & & & \\ 
0 & & 0 & 0 & $-\frac{1}{2}$ & 0 \\
1 & & 0 & $\frac{5}{24}$ & 0 & $-\frac{15}{128}$ \\
2 & &  $-\frac{1}{6}$ & 0 & $\frac{19}{192}$ & 0 \\
3 & & 0 & $-\frac{3}{32}$ & 0 &  \\
4 & & $\frac{1}{8}$ & 0 & $-\frac{7}{128}$ &  \\
5 & & 0 & $\frac{15}{256}$ & 0 &  \\
6 & & $-\frac{11}{96}$ & 0 &  &  \\
7 & &   0  &  &  &  \\
& &  & & &  \\
\hline
\end{tabular}
\caption{Annulus orbifold Gromov--Witten invariants $\langle
  \mathbf{1}_\frac{1}{4}^m   \mathbf{1}_\frac{1}{2}^n \rangle_0^{(3,1)}$of
  $[\bbC^3/\bbZ_4]$ in the asymmetric case at winding number $(3,1)$.}
\label{tab:annulus3}
\end{table}

\subsection{$B$-model, symmetric case}
\subsubsection{Disc and annulus invariants}
The $B$-model setup for the symmetric case is obtained via a combined $S$ and
$T$ transformation of the curve \eqref{eq:mirrorcurvez4}. A simple form is
obtained for framing $f=0$, where the derivative of the symmetric disc
potential is obtained from \eqref{eq:Adiscpotasymm} by sending $U\to 1/U$.
%
The Hori-Vafa differential now reads
\bea
d\lambda^{(\sigma),f=0}_{\bbC^3/\bbZ_4}(a_{\frac{1}{4}},a_{\frac{1}{2}},U)&=&\log
\Bigg(\frac{1}{2} \Bigg(U^2+a_{\frac{1}{4}} U+a_{\frac{1}{2}} \nn \\ &-&
\sqrt{\left(U^2+a_{\frac{1}{4}}
  U+a_{\frac{1}{2}}\right)^2-4}\Bigg)\Bigg)\frac{dU}{U} \nn \\
\label{eq:Adiscpotsymm}
\eea
and the open mirror map is again trivial
\beq
w=U
\eeq
Upon expanding in $w$ and plugging in the closed mirror map, we find
\beq
\begin{aligned}
\cF_{0,1}^{[\mathbb{C}^3/\bbZ_4], (\sigma)}(\tau_i,\tau_{\frac{1}{2}},w)  :=  \int^w
d\lambda^{(\sigma),f=0}_{\bbC^3/\bbZ_4}(a_{\frac{1}{4}}(\tau_i,
\tau_{\frac{1}{2}}),a_{\frac{1}{2}}(\tau_{\frac{1}{2}}),x) d x  \\ 
= i \Bigg[ \left(\left(\frac{1}{2}+\frac{3 \tau_{\frac{1}{2}}^2}{64}\right) \tau_{\frac{1}{4}}+\left(\frac{\tau_{\frac{1}{2}}}{384}+\frac{7
   \tau_{\frac{1}{2}}^3}{12288}\right)
\tau_{\frac{1}{4}}^3+O\left(\tau_{\frac{1}{4}}^4\right)\right)w  \\
 + 
   \left(\left(\frac{1}{4}+\frac{\tau_{\frac{1}{2}}^2}{32}\right)+\left(\frac{\tau_{\frac{1}{2}}}{32}-\frac{5
   \tau_{\frac{1}{2}}^3}{1536}\right)
   \tau_{\frac{1}{4}}^2+O\left(\tau_{\frac{1}{4}}^4\right) \right)w^2+O\left(w^3\right) \Bigg]  \\
\label{eq:bdiscsymm}
\end{aligned}
\eeq
The $B$-model potential thus has an expansion in rational numbers only up to a
$\pi/2$ phase. It would be interesting to track its origin in detail. \\
Likewise, the computation of annulus invariants requires basically no new ingredients
with respect to the asymmetric case, the only difference being that we have to
replace $w_1 \to 1/w_1$, $w_2 \to 1/w_2$ in the expression for the Bergmann
kernel. Conjecture \ref{conj:BKMP2} then allows us to compute
\beq
\begin{aligned}
\cF_{0,2}^{[\mathbb{C}^3/\bbZ_4], (\sigma)}(\tau_i,\tau_{\frac{1}{2}},w_1, w_2)  =   
\l(\frac{3
  \tau_{\frac{1}{4}}^2}{64}+\frac{\tau_{\frac{1}{2}}}{16}+\frac{3 \tau_{\frac{1}{2}}^2
  \tau_{\frac{1}{4}}^2}{128}+\frac{11 \tau_{\frac{1}{2}}^3}{768} +\dots\r) w_1
w_2  \\
+
\l(\frac{\tau_{\frac{1}{4}}}{16} +\frac{19 \tau_{\frac{1}{2}}^2
  \tau_{\frac{1}{4}}}{512}+\frac{43 \tau_{\frac{1}{2}}
  \tau_{\frac{1}{4}}^3}{3072}+ \frac{295 \tau_{\frac{1}{2}}^3 \tau_{\frac{1}{4}}^3}{32768}+\dots\r)(w_1^2 w_2 + w_2 w_1^2)   \\ + \dots
\label{eq:bannsymm}
\end{aligned}
\eeq
\subsection{$A$-model, symmetric case}

\label{sec:c3z4s}
In the symmetric case, we apply  formula (\ref{locformula}) to compute disc and annulus invariants for the orbifold $[\IC^3/\IZ_4]$, with $\alpha_0=2$, $\alpha_1=\alpha_2=1$, so the action is ineffective (with a $\IZ_2$ isotropy group) along the axis that gets doubled to become the zero section of the orbi-bundle. The torus weights
$$
(s_0,s_1,s_2)=\left(\frac{1}{2}, -a, a-\frac{1}{2}\right),
$$
can be specialized to $a=1/4$ to obtain symmetric weights in the fiber
directions. Once again, insertions that give rise to non-equivariant invariants correspond to the two age one orbifold cohomology classes, $\mathbf{1}_{\frac{1}{4}}, \mathbf{1}_{\frac{2}{4}}$.

The compatibility condition between degree and twisting (\ref{degtwist}) is $k\equiv d$ (mod $2$), and the disc function is:
\bea
\mathcal{D}^{(\sigma)}(d,1/4) &=& \frac{1}{2} \frac{ (\mathbf{1}_{\frac{d}{4}})^\vee}{\left\lfloor \frac{d}{2}\right\rfloor!}\left(\frac{\hbar}{d}\right)^{\mathrm{age}(\mathbf{1}_{\frac{d}{4}})}\frac{\Gamma(\frac{d}{2}+\langle \frac{d}{4}\rangle-\frac{d}{4})}{\Gamma(1-\langle \frac{d}{4}\rangle-\frac{d}{4})}
\nn \\ &+& \frac{1}{2}
\frac{ (\mathbf{1}_{\frac{d+2}{4}})^\vee}{\left\lfloor
  \frac{d}{2}\right\rfloor!}\left(\frac{\hbar}{d}\right)^{\mathrm{age}(\mathbf{1}_{\frac{d+2}{4}})}\frac{\Gamma(\frac{d}{2}+\langle
  \frac{d+2}{4}\rangle-\frac{d}{4})}{\Gamma(1-\langle
  \frac{d+2}{4}\rangle-\frac{d}{4})} \nn \\
\eea

\subsubsection{Disc Invariants}
We can then compute the potential for open disc invariants 


\begin{align}
\mathcal{F}_{0,1}^{[\mathbb{C}^3/\bbZ_4], (\sigma)}(\tau_\frac{1}{4},\tau_\frac{2}{4},w,\frac{1}{4}):=\sum_{m_1,m_2,d}\langle\mathbf{1}_{\frac{1}{4}}^{m_1}\mathbf{1}_{\frac{2}{4}}^{m_2};\frac{1}{4}\rangle^{d}_{0}\frac{\tau_{\frac{1}{4}}^{m_1}}{{m_1}!}\frac{\tau_{\frac{2}{4}}^{m_2}}{m_2!}\frac{w^d}{d!}= \nonumber\\
\sum_d\left[J\left(\frac{1}{2},-\frac{1}{4},-\frac{1}{4};0,\tau_\frac{1}{4},\tau_\frac{2}{4};\frac{1}{d}\right)\mathcal{D}^{(\sigma)}\l(d,\frac{1}{4}\r) \right]\frac{y^d}{d!}
\end{align}
where once again $J(s_0,s_1,s_2;\tau_{\frac{1}{4}},\tau_{\frac{2}{4}};z)$ is
the small $J$-function of the closed theory.
 
Explicit values for $n$-pointed disc invariants are shown and compared with the physical
predictions in Table \ref{tab:discsymm1}-\ref{tab:discsymm2}. The final result
agrees with the $B$-model prediction, apart from the usual sign ambiguity.

\begin{table}[h]
\centering
\begin{tabular}{|c|ccccc|}
\hline
 & $m$ & 1 & 3 & 5 & 7 \\  
\hline
$n$ & &  & & & \\ 
0 & &  $\frac{1}{2}$ & 0 & -$\frac{3}{128}$ & 0 \\
1 & &  0 & $\frac{1}{64}$ & 0 & $\frac{207}{4096}$ \\
2 & &  $\frac{3}{32}$ & 0 & -$\frac{49}{2048}$ & 0 \\
3 & &  0 & $\frac{21}{1024}$ & 0 & $\frac{12447}{65536}$ \\
4 & &  -$\frac{91}{512}$ & 0 & -$\frac{2247}{32768}$ & 0 \\
5 & &  0 & $\frac{361}{16384}$ & 0 & $\frac{1272327}{1048576}$ \\
6 & &  $\frac{1703}{8192}$ & 0 & -$\frac{191349}{524288}$ & 0 \\
7 & &  0 & $\frac{37661}{262144}$ & 0 & $\frac{202603527}{16777216}$ \\
& &  & & &  \\
\hline
\end{tabular}
\caption{Disc orbifold Gromov--Witten invariants $\langle
  \mathbf{1}_\frac{1}{4}^m   \mathbf{1}_\frac{1}{2}^n \rangle_0^{(1)}$of
  $[\bbC^3/\bbZ_4]$ in the symmetric case at winding number $1$.}
\label{tab:discsymm1}
\end{table}

\begin{table}[h]
\centering
\begin{tabular}{|c|ccccc|}
\hline
 & $m$ & 0 & 2 & 4 & 6 \\  
\hline
$n$ & &  & & & \\ 
0 & & $\frac{1}{4}$ & 0 & 0 & 0 \\
1 & & 0 & $\frac{1}{16}$ & 0 & -$\frac{9}{512}$ \\
2 & & $\frac{1}{16}$ & 0 & $\frac{1}{64}$ & 0 \\
3 & & 0 & $\frac{13}{128}$ & 0 & -$\frac{141}{2048}$ \\
4 & & $\frac{5}{64}$ & 0 & $\frac{19}{256}$ & 0 \\
5 & & 0 & $\frac{3}{8}$ & 0 & -$\frac{4189}{8192}$ \\
6 & & $\frac{61}{256}$ & 0 & $\frac{1137}{2048}$ & 0 \\
7 & & 0 & $\frac{10111}{4096}$ & 0 & -$\frac{98013}{16384}$ \\
& &  & & &  \\
\hline
\end{tabular}
\caption{Disc orbifold Gromov--Witten invariants $\langle
  \mathbf{1}_\frac{1}{4}^m   \mathbf{1}_\frac{1}{2}^n \rangle_0^{(2)}$of
  $[\bbC^3/\bbZ_4]$ in the symmetric case at winding number $2$.}
\label{tab:discsymm2}
\end{table}

\subsubsection{Annulus Invariants}

Our localization formula expresses annulus invariants in terms of the disc function and of compact invariants with two descendant insertions:

\begin{equation}
\begin{array}{l}
\langle \mathbf{1}^{m_0} \mathbf{1}_{\frac{1}{4}}^{m_1} \mathbf{1}_{\frac{2}{4}}^{m_{2}}\rangle_0^{d_1,d_2}=  \\ \\
\left(\frac{4}{\hbar}\right)^2\sum_{k_j \equiv d_j (\mbox{mod}\ 2)} \prod_{j=1}^2 \left(\varphi(k_j) \DD^{(s)}_{k_j}\left(d_j, 1/4\right)\right)\\
\int_{\mathcal{M}}\frac{e^{\mathrm{eq}}(\IE^\vee_{2}(1/2)\oplus\IE^\vee_{1}(-1/4)\oplus\IE^\vee_{1}(-1/4))}{ \left({\frac{\hbar}{d_1}}-\psi_1\right)\left({\frac{\hbar}{d_2}}-\psi_2\right)},
\end{array} 
\end{equation}
where $$
\mathcal{M}=\overline{\mathcal{M}}_{0,\sum m_j+2}(B\IZ_4,0;\mathbf{1}^{m_0} \mathbf{1}_{\frac{1}{4}}^{m_1} \mathbf{1}_{\frac{2}{4}}^{m_{2}},\mathbf{1}_{\frac{4-k_1}{4}}, \mathbf{1}_{\frac{4-k_2}{4}})
$$
and
$$
\varphi(k):=
\left\{
\begin{array}{cl}
\frac{\hbar^3}{32} & k=0 \\
1                  & k=1,3 \\
\frac{\hbar}{2}    &k=2
\end{array}
\right.
$$

As for the asymmetric case, we remove inductively the second descendant insertion using the genus $0$ topological recursion relations 
 and the string equation. In Table
\ref{tab:annsymm1}-\ref{tab:annsymm2} we collect the first values for these
invariants. Again, to the extent we have checked, we find agreement with the
$B$-model prediction.  

\begin{table}[h]
\centering
\begin{tabular}{|c|ccccc|}
\hline
 & $m$ & 1 & 3 & 5  & 7 \\  
\hline
$n$ & &  & & & \\ 
0 & &  0 & $\frac{3}{16}$ & 0 & -$\frac{21}{256}$  \\
1 & &  $\frac{1}{8}$ & 0 & $\frac{9}{256}$ & 0  \\
2 & &  0 & $\frac{3}{16}$ & 0 & -$\frac{159}{1024}$  \\
3 & &  $\frac{11}{64}$ & 0 & $\frac{63}{512}$ & 0  \\
4 & &  0 & $\frac{317}{512}$ & 0 &   \\
5 & &  $\frac{71}{128}$ & 0 & $\frac{28877}{65536}$ &   \\
6 & &  0 & -$\frac{4179}{8192}$ &  &   \\
7 & &  $\frac{6687}{2048}$ &  &  &  \\
& &  & & &  \\
\hline
\end{tabular}
\caption{Annulus orbifold Gromov--Witten invariants $\langle
  \mathbf{1}_\frac{1}{4}^m   \mathbf{1}_\frac{1}{2}^n \rangle_0^{(2)}$of
  $[\bbC^3/\bbZ_4]$ in the symmetric case at winding number $(1,1)$.}
\label{tab:annsymm1}
\end{table}

\begin{table}[h]
\centering
\begin{tabular}{|c|ccccc|}
\hline
 & $m$ & 0 & 2 & 4 & 6 \\  
\hline
$n$ & &  & & & \\ 
0 & & $\frac{1}{8}$ & 0 & -$\frac{3}{512}$ & 0 \\
1 & & 0 & $\frac{43}{256}$ & 0 & -$\frac{2859}{16384}$ \\
2 & & $\frac{19}{128}$ & 0 & $\frac{863}{8192}$ & 0 \\
3 & &  0 & $\frac{2655}{4096}$ & 0 &  \\
4 & & $\frac{1109}{2048}$ & 0 & $\frac{129513}{131072}$ &  \\
5 & & 0 & $\frac{306883}{65536}$ & 0 &  \\
6 & & $\frac{119719}{32768}$ & 0 &  &  \\
7 & &  0 &  &  &  \\
& &  & & &  \\
\hline
\end{tabular}
\caption{Annulus orbifold Gromov--Witten invariants $\langle
  \mathbf{1}_\frac{1}{4}^m   \mathbf{1}_\frac{1}{2}^n \rangle_0^{(2)}$of
  $[\bbC^3/\bbZ_4]$ in the symmetric case at winding number $(2,1)$.}
\label{tab:annsymm2}
\end{table}

\clearpage
\appendix
\section{The Eynard-Orantin recursion in the elliptic case}
\label{sec:appmodular}
We review some details of the Eynard-Orantin recursion specialized to the case when the support $\Gamma$ of the spectral curve $\mathcal{S}$ is a complex 2-torus, $\mathfrak{g}=1$, and $V$ realizes it as a degree 2 branched covering of $\bbC \bbP^1$. The Hori-Vafa differential \eqref{eq:horivafadiff} reads
\beq
\label{differ22}
d\lambda_X(U)= \log\left(\frac{P_2(U) \pm Y(U)}{2} \right)\frac{dU}{U}
\eeq
where
\beq
\hbox{deg }P_2(U)=2, \qquad Y(U)=\sqrt{P_2^2(U)-4}=\sqrt{U-q_1)(U-q_2)(U-q_3)(U-q_4)}
\label{defY}
\eeq
We have first of all that
\beq
\label{thxh}
d\lambda(U)-d\lambda(\bar U)=2 M(U) Y(U) d U
\eeq
where the so-called ``moment function'' $M(U)$ is given, after using the fact that $\log({P+Y})-\log({P-Y})= 2 \tanh^{-1}{(Y/P)}$, as
\beq
M(U)={1\over U Y(U)} \tanh^{-1} \biggl[ { {Y(U)} \over P_2(U)} \biggr],
\label{defM}
\eeq
Moreover, the one form $dE(p,q)$ can be written as \cite{Bouchard:2007ys}
\beq
\label{CjLambda}
d E_{W}(U)={1\over 2} {Y(W)\over Y(U)}\,\left(
{1\over U-W}-L C(W)
\right)\, d U
\eeq
where
\beq\label{defCj}
C(W):={1\over 2\pi i }\,\oint_{A} {d U \over {Y(U)}}\,{1\over U-W}, \qquad L^{-1}:= {1\over 2\pi i }\,\oint_{A} {d U \over {Y(U)}}
\eeq
We have assumed here that $W$ stays outside the contour $A$; when $W$ lies inside the contour $A$, $C(W)$ in (\ref{CjLambda}) should be replaced by its regularized version
\beq
\label{regc}
C^{\rm reg}(W)=C(W)-\frac{1}{Y(W)}
\eeq
%
%
%
%
Since $\Gamma_X$ is elliptic, it is possible to find closed form expressions for $C(U)$, $C_{reg}(U)$, $B(U,W)$ and $L$. We have 
\beq
\label{eq:C}
C(U)= {2 (q_2-q_3)\over \pi (U-q_3) (U-q_2) {\sqrt{(q_1-q_3)(q_2-q_4)}}} 
\biggl[  \Pi(n_4, k) + \frac{U-q_2}{q_2-q_3} K(k)\biggr]
\eeq
\beq
\label{eq:Creg}
C^{\rm reg}(U)=  {2 (q_3-q_2)\over \pi (U-q_3) (U-q_2) {\sqrt{(q_1-q_3)(q_2-q_4)}}} 
\biggl[  \Pi(n_1, k) + \frac{U-q_3}{q_3-q_2} K(k)\biggr] 
\eeq
\beq L^{-1} = \frac{2}{\sqrt{(q_1-q_3)(q_2-q_4)}}K\left[\frac{(q_1-q_2)(q_3-q_4)}{(q_1-q_3)(q_2-q_4)}\right]  \label{defL}\eeq
\bea B(U,W)&=&\frac{1}{Y(U)}\left[\frac{Y^2(U)}{2 Y(W) (U-W)^2}+\frac{(Y^2)'(U)}{4 Y(W) (W-U)}+\frac{A(U)}{4 Y(W)} \right]\nonumber \\ &+& \frac{1}{2(U-W)^2}
\label{eq:bergmann}
\eea
where 
\beq
k ={(q_1-q_2)(q_3-q_4)\over(q_1-q_3)(q_2-q_4)}, \qquad n_4= {(q_2-q_1)(U-q_3) \over (q_3-q_1)(U-q_2)}, \qquad n_1= {(q_4-q_3)(U-q_2) \over (q_4-q_2)(U-q_3)},
\label{ellmodulus}
\eeq
\beq A(U)=(U-q_1)(U-q_2)+(U-q_3)(U-q_4) +(q_1-q_3)(q_2-q_4)\frac{E(k)}{K(k)} \label{defA} \eeq
and $K(k)$, $E(k)$ and $\Pi(n, k)$ are the complete elliptic integrals of the first, second and third kind respectively.
%
%
\\

\noindent With these ingredients one can compute the residues in (\ref{eq:BKMPrecursion1})-\eqref{eq:BKMPrecursion3}. Given that $d E_{q}(p)/(d\lambda(q)-d\lambda(\bar q))$, as a function of $q$, 
is regular at the branch-points, all residues appearing in (\ref{eq:BKMPrecursion3}) will be linear combinations of the following {\it kernel differentials}
\bea
\label{kd}
\chi_i^{(n)}(p)&=&{\rm Res}_{q=q_i} \biggl(  \frac{d E_q (p)}{d\lambda(q)-d\lambda(\bar q)} {1\over (q-q_i)^n} \biggr) \nonumber \\
&=&{1\over  (n-1)!} {1\over  {Y(p)}} {d^{n-1} \over d q^{n-1}} \Biggl[{1\over 2 M(q) }\biggl( {1\over p-q} -L C(q)\biggl) \Biggr]_{q=q_i}
\eea
In (\ref{kd}), $C(p)$ should be replaced by $C_{\rm reg}(p)$ when $i=1,2$. 
%
\\

It is instructive to see the appearance in general of the propagators $G(q_1,q_2,q_3,q_4)$ as defined in \eqref{eq:Gprop}. Let $f:\bbC\to\bbC$ be a complex valued function with meromorphic square $f^2(x)$, and denote with $f^{(n)}_i$ the $(n+1)$-th coefficient in a Laurent expansion of $f(x)$ around $q_i$
\beq
f(x) = \sum_{n=-N_i}^\infty \frac{f_i^{(n+N_i)}}{(p-q_i)^{n/2}}
\eeq
Then Eq. \eqref{eq:BKMPrecursion1}-\eqref{eq:BKMPrecursion3} and (\ref{kd}) imply that the correlators will be a polynomial in the following four basic building blocks
\beq
\label{objrec}
M_i^{(n)}, \qquad  A^{(n)}_i, \qquad \left(\frac{1}{Y}\right)_i^{(n)}, \qquad \mathcal{C}^{(n)}_i
\eeq
where, we  have defined
\beq
\mathcal{C}^{(n)}_i=\left\{\begin{array}{cc} C^{(n)}_{{\rm reg}, i} & \mathrm{for }\quad  i=1,2 \\ C^{(n)}_i & \mathrm{for } \quad i=3,4\end{array}\right.\eeq
It is immediate to see that the residue computation involving $M_i^{(n)}$ will always yield an algebraic function of the ``bare'' complex moduli, that is the coefficients of $P_2$. This means that they have degree zero as a polynomial in $G(\{q_i\})$. On the other hand, they are the only ones who bring a dependence on the marked functions of the spectral curve $\mathcal{S}_X$: all the others only depend on differences of branch points $q_i$, which (perhaps up to a rescaling of $U$ and $V$) leads to functions of the elliptic modulus of $\Gamma_X$ which are linear in $G(\{q_i\})$. This is apparent for  $A^{(n)}_i$ and $(1/Y)_i^{(n)}$ from formulae (\ref{defY}) and (\ref{defA}), while the case of $\mathcal{C}^{(n)}_i$ follows from the fact that
$$ \partial_x \Pi(x,y)= \frac{x E(y)+(y-x) K(y)+\left(x^2-y\right) \Pi (x,y)}{2
   (x-1) x (y-x)}
$$
The above formula implies that 
\beq
\partial^{(n)}_x \Pi(x,y)=A_n(x,y) K(y) +  B_n(x,y) E(y)+ C_n(x,y) \Pi(x,y)
\label{eq:derpi}
\eeq
where $A_n$, $B_n$ and $C_n$ are rational functions of $x$ and $y$. From (\ref{ellmodulus}), to compute $\mathcal{C}^{(n)}_i$, we need to evaluate these expressions when $n_1$ (resp. $n_4$) equals either 0 or $k$. But using
\beq
 \Pi(0,y)=K(y), \qquad \Pi(y,y)=\frac{E(y)}{1-y}
\eeq
we conclude that 
\beq
\mathcal{C}^{(n)}_i=R^{(n)}_1\l(\{q_i\}\r) K(k) + R^{(n)}_2\l(\{q_i\}\r) E(k)
\eeq
for two sequences of rational functions $R^{(n)}_i$. Notice that by (\ref{kd}), $\mathcal{C}^{(n)}_i$ always appears multiplied by $L$ in the recursion; therefore, from (\ref{defL}), $L \mathcal{C}^{(n)}_i$ is linear in $G(\{q_i\})$.

\clearpage
\section{Mirror symmetry predictions of open orbifold $GW$ invariants of
  $[\bbC^3/\bbZ_4]$ in the asymmetric case}
\label{sec:numbers}
\begin{table}[h]
\centering
\begin{tabular}{|c|ccccc|}
\hline
 & $m$ & 1 & 3 & 5 & 7  \\ 
\hline
$n$ & &  & & &  \\ 
0 & & 0 & -$\frac{3}{64}$ & 0 & $\frac{189}{4096}$  \\ 
1 & &  $\frac{1}{32}$ & 0 & -$\frac{33}{2048}$ & 0  \\ 
2 & &  0 & $\frac{11}{1024}$ & 0 & $\frac{14547}{65536}$  \\ 
3 & &  -$\frac{7}{512}$ & 0 & -$\frac{1989}{32768}$ & 0  \\ 
4 & &  0 & $\frac{353}{16384}$ & 0 & $\frac{1809801}{1048576}$  \\ 
5 & &  -$\frac{79}{8192}$ & 0 & -$\frac{218993}{524288}$ & 0  \\ 
6 & &  0 & $\frac{33711}{262144}$ & 0 & $\frac{330787647}{16777216}$  \\ 
7 & &  -$\frac{7287}{131072}$ & 0 & -$\frac{36190149}{8388608}$ & 0  \\ 
8 & &  0 & $\frac{4907493}{4194304}$ & 0 & $\frac{84814988181}{268435456}$  \\ 
9 & &  -$\frac{889439}{2097152}$ & 0 & -$\frac{8528369313}{134217728}$ & 0  \\ 
10 & &  0 & $\frac{1045989811}{67108864}$ & 0 & $\frac{29188217357547}{4294967296}$  \\ 
11 & &  -$\frac{167510567}{33554432}$ & 0 & -$\frac{2728134070309}{2147483648}$ & 0  \\ 
12 & &  0 & $\frac{307481197833}{1073741824}$ & 0 & $\frac{13004327932052961}{68719476736}$  \\ 
& &  & & &  \\ 
\hline
\end{tabular}
\caption{Predictions for $g=0$, $h=3$ open orbifold Gromov--Witten invariants of $[\bbC^3/\bbZ_4]$ at winding number $(1,1,1)$.}
\end{table}

\begin{table}[h]
\centering
\begin{tabular}{|c|cccccc|}
\hline
 & $m$ & 0 & 2 & 4 & 6 & 8  \\ 
\hline
$n$ & &  & & & &  \\ 
0 & & $\frac{1}{2}$ & 0 & $\frac{9}{64}$ & 0 & -$\frac{21}{128}$  \\ 
1 & &  0 & -$\frac{1}{16}$ & 0 & $\frac{3}{64}$ & 0  \\ 
2 & &  $\frac{1}{16}$ & 0 & -$\frac{7}{256}$ & 0 & -$\frac{4377}{8192}$  \\ 
3 & &  0 & $\frac{3}{128}$ & 0 & $\frac{87}{1024}$ & 0  \\ 
4 & &  -$\frac{1}{32}$ & 0 & -$\frac{3}{512}$ & 0 & -$\frac{64593}{16384}$  \\ 
5 & &  0 & -$\frac{7}{512}$ & 0 & $\frac{4509}{8192}$ & 0  \\ 
6 & &  $\frac{17}{512}$ & 0 & -$\frac{111}{2048}$ & 0 & -$\frac{5600217}{131072}$  \\ 
7 & &  0 & -$\frac{1}{64}$ & 0 & $\frac{170919}{32768}$ & 0  \\ 
8 & &  $\frac{9}{512}$ & 0 & -$\frac{11867}{32768}$ & 0 & -$\frac{85512669}{131072}$  \\ 
9 & &  0 & -$\frac{423}{2048}$ & 0 & $\frac{4679277}{65536}$ & 0  \\ 
10 & &  $\frac{1091}{4096}$ & 0 & -$\frac{493299}{131072}$ & 0 & -$\frac{14083706541}{1048576}$  \\ 
11 & &  0 & -$\frac{176659}{65536}$ & 0 & $\frac{701236689}{524288}$ & 0  \\ 
12 & &  $\frac{22219}{8192}$ & 0 & -$\frac{7019643}{131072}$ & 0 & -$\frac{753800096679}{2097152}$  \\ 
& &  & & & &  \\ 
\hline
\end{tabular}
\caption{Predictions for $g=0$, $h=3$ open orbifold Gromov--Witten invariants of $[\bbC^3/\bbZ_4]$ at winding number $(2,1,1)$.}
\end{table}

\begin{table}[h]
\centering
\begin{tabular}{|c|ccccc|}
\hline
 & $m$ & 1 & 3 & 5 & 7  \\ 
\hline
$n$ & &  & & &  \\ 
0 & &  -$\frac{3}{4}$ & 0 & -$\frac{111}{256}$ & 0  \\ 
1 & &  0 & $\frac{15}{128}$ & 0 & -$\frac{2595}{8192}$  \\ 
2 & &  -$\frac{1}{64}$ & 0 & $\frac{603}{4096}$ & 0  \\ 
3 & &  0 & -$\frac{177}{2048}$ & 0 & -$\frac{40479}{131072}$  \\ 
4 & &  $\frac{65}{1024}$ & 0 & $\frac{1245}{65536}$ & 0  \\ 
5 & &  0 & $\frac{895}{32768}$ & 0 & -$\frac{4055235}{2097152}$  \\ 
6 & &  -$\frac{541}{16384}$ & 0 & $\frac{343623}{1048576}$ & 0  \\ 
7 & &  0 & -$\frac{68737}{524288}$ & 0 & -$\frac{568999599}{33554432}$  \\ 
8 & &  $\frac{34245}{262144}$ & 0 & $\frac{39364105}{16777216}$ & 0  \\ 
9 & &  0 & -$\frac{5665425}{8388608}$ & 0 & -$\frac{119917956675}{536870912}$  \\ 
10 & &  $\frac{1984519}{4194304}$ & 0 & $\frac{7678005843}{268435456}$ & 0  \\ 
11 & &  0 & -$\frac{1112041297}{134217728}$ & 0 & -$\frac{35129545858719}{8589934592}$  \\ 
12 & &  $\frac{386924425}{67108864}$ & 0 & $\frac{2109027490965}{4294967296}$ & 0  \\ 
& &  & & &  \\ 
\hline
\end{tabular}
\caption{Predictions for $g=0$, $h=3$ open orbifold Gromov--Witten invariants of $[\bbC^3/\bbZ_4]$ at winding number $(2,2,1)$.}
\end{table}

\begin{table}[h]
\centering
\begin{tabular}{|c|ccccc|}
\hline
 & $m$ & 1 & 3 & 5 & 7 \\ 
\hline
$n$ & &  & & &  \\ 
0 & & -$\frac{2}{3}$ & 0 & -$\frac{19}{32}$ & 0  \\ 
1 & &  0 & $\frac{35}{192}$ & 0 & -$\frac{1655}{4096}$  \\ 
2 & &  -$\frac{1}{16}$ & 0 & $\frac{599}{3072}$ & 0  \\ 
3 & &  0 & -$\frac{127}{1024}$ & 0 & -$\frac{34639}{65536}$  \\ 
4 & &  $\frac{7}{64}$ & 0 & $\frac{71}{1024}$ & 0  \\ 
5 & &  0 & $\frac{385}{16384}$ & 0 & -$\frac{4900175}{1048576}$  \\ 
6 & &  -$\frac{733}{12288}$ & 0 & $\frac{263963}{262144}$ & 0  \\ 
7 & &  0 & -$\frac{280541}{786432}$ & 0 & -$\frac{972023479}{16777216}$  \\ 
8 & &  $\frac{8599}{32768}$ & 0 & $\frac{71107289}{6291456}$ & 0  \\ 
9 & &  0 & -$\frac{12976415}{4194304}$ & 0 & -$\frac{276764938375}{268435456}$  \\ 
10 & &  $\frac{1372579}{1048576}$ & 0 & $\frac{12755007193}{67108864}$ & 0  \\ 
11 & &  0 & -$\frac{3203261567}{67108864}$ & 0 & -$\frac{105795705480319}{4294967296}$  \\ 
12 & &  $\frac{28966097}{1572864}$ & 0 & $\frac{1153922108479}{268435456}$ & 0  \\ 
& &  & & &  \\ 
\hline
\end{tabular}
\caption{Predictions for $g=0$, $h=3$ open orbifold Gromov--Witten invariants of $[\bbC^3/\bbZ_4]$ at winding number $(3,1,1)$.}
\end{table}

\newpage
\pagebreak
\newpage


\begin{table}[h]
\centering
\begin{tabular}{|c|ccccc|}
\hline
 & $m$ & 1 & 3 & 5 & 7  \\ 
\hline
$n$ & &  & & &  \\ 
0 & & $\frac{1}{48}$ & 0 & $\frac{9}{1024}$ & 0  \\ 
1 & & 0 & -$\frac{5}{1536}$ & 0 & -$\frac{3375}{32768}$  \\ 
2 & & $\frac{1}{768}$ & 0 & $\frac{1367}{49152}$ & 0  \\ 
3 & & 0 & -$\frac{81}{8192}$ & 0 & -$\frac{326497}{524288}$  \\ 
4 & & $\frac{65}{12288}$ & 0 & $\frac{40345}{262144}$ & 0  \\ 
5 & & 0 & -$\frac{19145}{393216}$ & 0 & -$\frac{50714835}{8388608}$  \\ 
6 & & $\frac{4321}{196608}$ & 0 & $\frac{5760669}{4194304}$ & 0  \\ 
7 & & 0 & -$\frac{2469623}{6291456}$ & 0 & -$\frac{11529490917}{134217728}$  \\ 
8 & & $\frac{490945}{3145728}$ & 0 & $\frac{3642090395}{201326592}$ & 0  \\ 
9 & & 0 & -$\frac{158835215}{33554432}$ & 0 & -$\frac{3604297162935}{2147483648}$  \\ 
10 & & $\frac{85184641}{50331648}$ & 0 & $\frac{354513303549}{1073741824}$ & 0  \\ 
11 & & 0 & -$\frac{128738647003}{1610612736}$ & 0 & -$\frac{1481653476327337}{34359738368}$  \\ 
12 & & $\frac{21004177025}{805306368}$ & 0 & $\frac{137005640391385}{17179869184}$ & 0  \\ 
& &  & & &  \\ 
\hline
\end{tabular}
\caption{Predictions for $g=1$, $h=1$ open orbifold Gromov--Witten invariants of $[\bbC^3/\bbZ_4]$ at winding number 1.}
\end{table}

\begin{table}[h]
\centering
\begin{tabular}{|c|cccccc|}
\hline
 & $m$ & 0 & 2 & 4 & 6 & 8  \\ 
\hline
$n$ & &  & & & &  \\ 
0 & & 0 & -$\frac{1}{32}$ & 0 & -$\frac{3}{512}$ & 0  \\ 
1 & &  $\frac{1}{48}$ & 0 & -$\frac{1}{512}$ & 0 & $\frac{19}{64}$  \\ 
2 & &  0 & $\frac{1}{192}$ & 0 & -$\frac{27}{512}$ & 0  \\ 
3 & &  -$\frac{5}{384}$ & 0 & $\frac{5}{512}$ & 0 & $\frac{121865}{65536}$  \\ 
4 & &  0 & -$\frac{5}{3072}$ & 0 & -$\frac{315}{1024}$ & 0  \\ 
5 & &  $\frac{1}{768}$ & 0 & $\frac{1277}{24576}$ & 0 & $\frac{1175231}{65536}$  \\ 
6 & &  0 & -$\frac{11}{2048}$ & 0 & -$\frac{89893}{32768}$ & 0  \\ 
7 & &  -$\frac{85}{12288}$ & 0 & $\frac{6995}{16384}$ & 0 & $\frac{131562305}{524288}$  \\ 
8 & &  0 & -$\frac{275}{6144}$ & 0 & -$\frac{9362985}{262144}$ & 0  \\ 
9 & &  -$\frac{389}{12288}$ & 0 & $\frac{1350073}{262144}$ & 0 & $\frac{631910777}{131072}$  \\ 
10 & &  0 & -$\frac{25459}{49152}$ & 0 & -$\frac{337003153}{524288}$ & 0  \\ 
11 & &  -$\frac{27155}{98304}$ & 0 & $\frac{4250765}{49152}$ & 0 & $\frac{2044307220305}{16777216}$  \\ 
12 & &  0 & -$\frac{4343505}{524288}$ & 0 & -$\frac{16051763495}{1048576}$ & 0  \\ 
& &  & & & &  \\ 
\hline
\end{tabular}
\caption{Predictions for $g=1$, $h=1$ open orbifold Gromov--Witten invariants of $[\bbC^3/\bbZ_4]$ at winding number $2$.}
\end{table}

\begin{table}[h]
\centering
\begin{tabular}{|c|ccccc|}
\hline
 & $m$ & 1 & 3 & 5 & 7  \\ 
\hline
$n$ & &  & & &  \\ 
0 & & 0 & $\frac{1}{12}$ & 0 & -$\frac{21}{256}$  \\ 
1 & &  -$\frac{5}{144}$ & 0 & $\frac{115}{3072}$ & 0  \\ 
2 & &  0 & -$\frac{59}{2304}$ & 0 & $\frac{9679}{49152}$  \\ 
3 & &  $\frac{53}{2304}$ & 0 & -$\frac{1903}{49152}$ & 0  \\ 
4 & &  0 & $\frac{305}{18432}$ & 0 & $\frac{219465}{131072}$  \\ 
5 & &  -$\frac{235}{12288}$ & 0 & -$\frac{716555}{2359296}$ & 0  \\ 
6 & &  0 & $\frac{4819}{65536}$ & 0 & $\frac{240631049}{12582912}$  \\ 
7 & &  -$\frac{16007}{589824}$ & 0 & -$\frac{14352681}{4194304}$ & 0  \\ 
8 & &  0 & $\frac{962165}{1179648}$ & 0 & $\frac{7715651635}{25165824}$  \\ 
9 & &  -$\frac{2979965}{9437184}$ & 0 & -$\frac{3575613975}{67108864}$ & 0  \\ 
10 & &  0 & $\frac{1821378401}{150994944}$ & 0 & $\frac{21402084232819}{3221225472}$  \\ 
11 & &  -$\frac{207837889}{50331648}$ & 0 & -$\frac{10779639149749}{9663676416}$ & 0  \\ 
12 & &  0 & $\frac{32205472535}{134217728}$ & 0 & $\frac{1609762782468295}{8589934592}$  \\ 
& &  & & &  \\ 
\hline
\end{tabular}
\caption{Predictions for $g=1$, $h=1$ open orbifold Gromov--Witten invariants of $[\bbC^3/\bbZ_4]$ at winding number 3.}
\end{table}

\begin{table}[h]
\centering
\begin{tabular}{|c|cccccc|}
\hline
 & $m$ & 0 & 2 & 4 & 6 & 8  \\ 
\hline
$n$ & &  & & & &  \\ 
0 & &  $\frac{1}{3}$ & 0 & -$\frac{5}{16}$ & 0 & $\frac{315}{256}$  \\ 
1 & &  0 & $\frac{3}{32}$ & 0 & -$\frac{181}{512}$ & 0  \\ 
2 & &  -$\frac{1}{24}$ & 0 & $\frac{39}{256}$ & 0 & -$\frac{109}{128}$  \\ 
3 & &  0 & -$\frac{11}{128}$ & 0 & $\frac{387}{2048}$ & 0  \\ 
4 & &  $\frac{7}{96}$ & 0 & -$\frac{3}{32}$ & 0 & -$\frac{162943}{16384}$  \\ 
5 & &  0 & $\frac{233}{3072}$ & 0 & $\frac{11689}{8192}$ & 0  \\ 
6 & &  -$\frac{31}{384}$ & 0 & -$\frac{403}{2048}$ & 0 & -$\frac{8297873}{65536}$  \\ 
7 & &  0 & -$\frac{331}{12288}$ & 0 & $\frac{286821}{16384}$ & 0  \\ 
8 & &  $\frac{187}{1536}$ & 0 & -$\frac{31927}{12288}$ & 0 & -$\frac{70090611}{32768}$  \\ 
9 & &  0 & $\frac{1389}{4096}$ & 0 & $\frac{73463763}{262144}$ & 0  \\ 
10 & &  $\frac{389}{6144}$ & 0 & -$\frac{5027977}{131072}$ & 0 & -$\frac{25016794729}{524288}$  \\ 
11 & &  0 & $\frac{193951}{49152}$ & 0 & $\frac{6224145569}{1048576}$ & 0  \\ 
12 & &  $\frac{47767}{24576}$ & 0 & -$\frac{99929913}{131072}$ & 0 & -$\frac{2869421365529}{2097152}$  \\ 
& &  & & & &  \\ 
\hline
\end{tabular}
\caption{Predictions for $g=1$, $h=1$ open orbifold Gromov--Witten invariants of $[\bbC^3/\bbZ_4]$ at winding number $4$.}
\end{table}

\newpage
\pagebreak
\newpage


\begin{table}[h]
\centering
\begin{tabular}{|c|ccccc|}
\hline
 & $m$ & 0 & 2 & 4 & 6  \\ 
\hline
$n$ & &  & & &  \\ 
0 & & 0 & $\frac{1}{192}$ & 0 & -$\frac{83}{2048}$  \\ 
1 & &  -$\frac{1}{96}$ & 0 & $\frac{19}{1536}$ & 0  \\ 
2 & &  0 & -$\frac{3}{512}$ & 0 & -$\frac{1313}{6144}$  \\ 
3 & &  $\frac{1}{128}$ & 0 & $\frac{59}{1024}$ & 0  \\ 
4 & &  0 & -$\frac{127}{6144}$ & 0 & -$\frac{116319}{65536}$  \\ 
5 & &  $\frac{11}{1024}$ & 0 & $\frac{10619}{24576}$ & 0  \\ 
6 & &  0 & -$\frac{3329}{24576}$ & 0 & -$\frac{1426777}{65536}$  \\ 
7 & &  $\frac{761}{12288}$ & 0 & $\frac{238597}{49152}$ & 0  \\ 
8 & &  0 & -$\frac{5593}{4096}$ & 0 & -$\frac{1171872737}{3145728}$  \\ 
9 & &  $\frac{4423}{8192}$ & 0 & $\frac{20167831}{262144}$ & 0  \\ 
10 & &  0 & -$\frac{15545773}{786432}$ & 0 & -$\frac{4454455581}{524288}$  \\ 
11 & &  $\frac{228811}{32768}$ & 0 & $\frac{2578302709}{1572864}$ & 0  \\ 
12 & &  0 & -$\frac{1225078949}{3145728}$ & 0 & -$\frac{4181947560489}{16777216}$  \\ 
& &  & & &  \\ 
\hline
\end{tabular}
\caption{Predictions for $g=1$, $h=2$ open orbifold Gromov--Witten invariants of $[\bbC^3/\bbZ_4]$ at winding number $(1,1)$.}
\end{table}

\begin{table}[h]
\centering
\begin{tabular}{|c|ccccc|}
\hline
 & $m$ & 1 & 3 & 5 & 7  \\ 
\hline
$n$ & &  & & &  \\ 
0 & & 0 & -$\frac{3}{128}$ & 0 & $\frac{749}{8192}$  \\ 
1 & & $\frac{1}{64}$ & 0 & -$\frac{73}{4096}$ & 0  \\ 
2 & & 0 & $\frac{11}{2048}$ & 0 & $\frac{59347}{131072}$  \\ 
3 & & -$\frac{13}{3072}$ & 0 & -$\frac{4277}{65536}$ & 0  \\ 
4 & & 0 & $\frac{323}{98304}$ & 0 & $\frac{7763577}{2097152}$  \\ 
5 & & $\frac{161}{16384}$ & 0 & -$\frac{1512059}{3145728}$ & 0  \\ 
6 & & 0 & $\frac{14031}{524288}$ & 0 & $\frac{1484285887}{33554432}$  \\ 
7 & & $\frac{28147}{786432}$ & 0 & -$\frac{87714437}{16777216}$ & 0  \\ 
8 & & 0 & $\frac{6282863}{25165824}$ & 0 & $\frac{395292407237}{536870912}$  \\ 
9 & & $\frac{4059523}{12582912}$ & 0 & -$\frac{21529218793}{268435456}$ & 0  \\ 
10 & & 0 & $\frac{1355800153}{402653184}$ & 0 & $\frac{140469708824427}{8589934592}$  \\ 
11 & & $\frac{277853969}{67108864}$ & 0 & -$\frac{21346688578591}{12884901888}$ & 0  \\ 
12 & & 0 & $\frac{132430149801}{2147483648}$ & 0 & $\frac{64325487060690897}{137438953472}$  \\ 
& &  & & &  \\ 
\hline
\end{tabular}
\caption{Predictions for $g=1$, $h=2$ open orbifold Gromov--Witten invariants of $[\bbC^3/\bbZ_4]$ at winding number $(2,1)$.}
\end{table}

\begin{table}[h]
\centering
\begin{tabular}{|c|cccccc|}
\hline
 & $m$ & 0 & 2 & 4 & 6 & 8  \\ 
\hline
$n$ & &  & & & &  \\ 
0 & & -$\frac{1}{4}$ & 0 & $\frac{7}{128}$ & 0 & -$\frac{63}{256}$  \\ 
1 & & 0 & -$\frac{1}{96}$ & 0 & $\frac{19}{512}$ & 0  \\ 
2 & & -$\frac{1}{32}$ & 0 & -$\frac{19}{1536}$ & 0 & -$\frac{15543}{16384}$  \\ 
3 & & 0 & $\frac{5}{768}$ & 0 & $\frac{85}{2048}$ & 0  \\ 
4 & & $\frac{1}{64}$ & 0 & $\frac{73}{3072}$ & 0 & -$\frac{254719}{32768}$  \\ 
5 & & 0 & -$\frac{47}{3072}$ & 0 & $\frac{2213}{6144}$ & 0  \\ 
6 & & -$\frac{17}{1024}$ & 0 & $\frac{391}{4096}$ & 0 & -$\frac{23761463}{262144}$  \\ 
7 & & 0 & -$\frac{25}{768}$ & 0 & $\frac{246515}{65536}$ & 0  \\ 
8 & & -$\frac{9}{1024}$ & 0 & $\frac{187937}{196608}$ & 0 & -$\frac{383755867}{262144}$  \\ 
9 & & 0 & -$\frac{4013}{12288}$ & 0 & $\frac{14299413}{262144}$ & 0  \\ 
10 & & -$\frac{1091}{8192}$ & 0 & $\frac{10594097}{786432}$ & 0 & -$\frac{66113744979}{2097152}$  \\ 
11 & & 0 & -$\frac{1594085}{393216}$ & 0 & $\frac{3320919305}{3145728}$ & 0  \\ 
12 & & -$\frac{22219}{16384}$ & 0 & $\frac{68613211}{262144}$ & 0 & -$\frac{3672961920137}{4194304}$  \\ 
& &  & & & &  \\ 
\hline
\end{tabular}
\caption{Predictions for $g=1$, $h=2$ open orbifold Gromov--Witten invariants of $[\bbC^3/\bbZ_4]$ at winding number $(2,2)$.}
\end{table}

\begin{table}[h]
\centering
\begin{tabular}{|c|ccccc|}
\hline
 & $m$ & 0 & 2 & 4 & 6  \\ 
\hline
$n$ & &  & & &  \\ 
0 & & -$\frac{4}{9}$ & 0 & $\frac{1}{8}$ & 0  \\ 
1 & &  0 & -$\frac{29}{576}$ & 0 & $\frac{189}{2048}$  \\ 
2 & &  $\frac{5}{144}$ & 0 & -$\frac{95}{2304}$ & 0  \\ 
3 & &  0 & $\frac{17}{512}$ & 0 & $\frac{6733}{24576}$  \\ 
4 & &  -$\frac{1}{18}$ & 0 & -$\frac{59}{1536}$ & 0  \\ 
5 & &  0 & -$\frac{13}{18432}$ & 0 & $\frac{215983}{65536}$  \\ 
6 & &  $\frac{115}{4608}$ & 0 & -$\frac{1335}{2048}$ & 0  \\ 
7 & &  0 & $\frac{16493}{73728}$ & 0 & $\frac{13351407}{262144}$  \\ 
8 & &  -$\frac{433}{2304}$ & 0 & -$\frac{179953}{18432}$ & 0  \\ 
9 & &  0 & $\frac{8527}{3072}$ & 0 & $\frac{1095436053}{1048576}$  \\ 
10 & &  -$\frac{50965}{36864}$ & 0 & -$\frac{25251695}{131072}$ & 0  \\ 
11 & &  0 & $\frac{117348941}{2359296}$ & 0 & $\frac{348308890001}{12582912}$  \\ 
12 & &  -$\frac{387979}{18432}$ & 0 & -$\frac{3841015873}{786432}$ & 0  \\ 
& &  & & &  \\ 
\hline
\end{tabular}
\caption{Predictions for $g=1$, $h=2$ open orbifold Gromov--Witten invariants of $[\bbC^3/\bbZ_4]$ at winding number $(3,1)$.}
\end{table}

\newpage
\pagebreak
\newpage

\begin{table}[h]
\centering
\begin{tabular}{|c|cccc|}
\hline
 & $m$ & 1 & 3 & 5  \\ 
\hline
$n$ & &  & &  \\ 
0 & &  $\frac{1}{3840}$ & 0 & -$\frac{851}{81920}$  \\ 
1 & &  0 & $\frac{391}{122880}$ & 0  \\ 
2 & &  -$\frac{91}{61440}$ & 0 & -$\frac{235817}{3932160}$  \\ 
3 & &  0 & $\frac{32689}{1966080}$ & 0  \\ 
4 & &  -$\frac{6263}{983040}$ & 0 & -$\frac{11465707}{20971520}$  \\ 
5 & &  0 & $\frac{1451137}{10485760}$ & 0  \\ 
6 & &  -$\frac{739891}{15728640}$ & 0 & -$\frac{7332916417}{1006632960}$  \\ 
7 & &  0 & $\frac{855699469}{503316480}$ & 0  \\ 
8 & &  -$\frac{394660109}{754974720}$ & 0 & -$\frac{2156215517801}{16106127360}$  \\ 
9 & &  0 & $\frac{233999737631}{8053063680}$ & 0  \\ 
10 & &  -$\frac{32994415691}{4026531840}$ & 0 & -$\frac{279022560888339}{85899345920}$  \\ 
11 & & 0 & $\frac{85088983138249}{128849018880}$ & 0  \\ 
12 & & -$\frac{3705845271181}{21474836480}$ & 0 &
-$\frac{414824442483351281}{4123168604160}$  \\ 
& &  & &  \\ 
\hline
\end{tabular}
\caption{Predictions for $g=2$, $h=1$ open orbifold Gromov--Witten invariants of $[\bbC^3/\bbZ_4]$ at winding number 1.}
\end{table}

\begin{table}[h]
\centering
\begin{tabular}{|c|ccccc|}
\hline
 & $m$ & 0 & 2 & 4 & 6  \\ 
\hline
$n$ & &  & & &  \\ 
0 & & 0 & -$\frac{23}{3840}$ & 0 & $\frac{397}{20480}$  \\ 
1 & & $\frac{37}{5760}$ & 0 & -$\frac{9}{2560}$ & 0  \\ 
2 & & 0 & $\frac{47}{46080}$ & 0 & $\frac{6763}{61440}$  \\ 
3 & & -$\frac{37}{23040}$ & 0 & -$\frac{1837}{122880}$ & 0  \\ 
4 & & 0 & $\frac{11}{73728}$ & 0 & $\frac{404921}{393216}$  \\ 
5 & & $\frac{599}{184320}$ & 0 & -$\frac{97087}{737280}$ & 0  \\ 
6 & & 0 & $\frac{2443}{491520}$ & 0 & $\frac{108494989}{7864320}$  \\ 
7 & & $\frac{10501}{737280}$ & 0 & -$\frac{2175833}{1310720}$ & 0  \\ 
8 & & 0 & $\frac{22249}{294912}$ & 0 & $\frac{1589644841}{6291456}$  \\ 
9 & & $\frac{191287}{1474560}$ & 0 & -$\frac{9401217}{327680}$ & 0  \\ 
10 & & 0 & $\frac{517601}{368640}$ & 0 & $\frac{15956291063}{2621440}$  \\ 
11 & & $\frac{10024913}{5898240}$ & 0 & -$\frac{123550208597}{188743680}$ & 0  \\ 
12 & & 0 & $\frac{414232039}{12582912}$ & 0 & $\frac{9407394255163}{50331648}$  \\ 
& &  & &  &  \\ 
\hline
\end{tabular}
\caption{Predictions for $g=2$, $h=1$ open orbifold Gromov--Witten invariants of $[\bbC^3/\bbZ_4]$ at winding number 2.}
\end{table}

\begin{table}[h]
\centering
\begin{tabular}{|c|cccc|}
\hline
 & $m$ & 1 & 3 & 5 \\ 
\hline
$n$ & &  & &  \\ 
0 & & 0 & $\frac{221}{5760}$ & 0  \\ 
1 & & -$\frac{643}{34560}$ & 0 & $\frac{19633}{737280}$  \\ 
2 & & 0 & -$\frac{7651}{552960}$ & 0  \\ 
3 & & $\frac{7231}{552960}$ & 0 & $\frac{907759}{11796480}$  \\ 
4 & & 0 & -$\frac{16787}{1105920}$ & 0  \\ 
5 & & $\frac{24737}{8847360}$ & 0 & $\frac{583154159}{566231040}$  \\ 
6 & & 0 & -$\frac{10955107}{47185920}$ & 0  \\ 
7 & & $\frac{12947971}{141557760}$ & 0 & $\frac{53326890019}{3019898880}$  \\ 
8 & & 0 & -$\frac{4233264583}{1132462080}$ & 0  \\ 
9 & & $\frac{2754790277}{2264924160}$ & 0 & $\frac{2105326035057}{5368709120}$  \\ 
10 & & 0 & -$\frac{2843659084991}{36238786560}$ & 0  \\ 
11 & & $\frac{825581791111}{36238786560}$ & 0 &
$\frac{25692970234497637}{2319282339840}$  \\ 
12 & & 0 & -$\frac{101535236275903}{48318382080}$ & 0  \\ 
& &  & &  \\ 
\hline
\end{tabular}
\caption{Predictions for $g=2$, $h=1$ open orbifold Gromov--Witten invariants of $[\bbC^3/\bbZ_4]$ at winding number 3.}
\label{tab:lastpred}
\end{table}

\clearpage

\bibliographystyle{alpha}
\bibliography{miabiblio}

\end{document}